\documentclass[11pt]{article}
\usepackage{latexsym, epsfig, amssymb, amsmath, amsthm, graphicx, mathrsfs}
\usepackage{anysize}
\usepackage[OT1]{fontenc}

\renewcommand{\baselinestretch}{1.2}
\marginsize{1.2in}{1in}{0.55in}{1.55in} %
\makeatletter
\def\singlespace{\def\baselinestretch{1}\@normalsize}

\newtheorem{assumption}{Condition}
\newtheorem{proposition}{Proposition}
\newtheorem{theorem}{Theorem}
\newtheorem{corollary}{Corollary}

\makeatother

\newcommand{\aphi}{\phi}

\newcommand{\ba}{\mbox{\bf a}}
\newcommand{\bb}{\mbox{\bf b}}
\newcommand{\be}{\mbox{\bf e}}

\newcommand{\br}{\mbox{\bf r}}
\newcommand{\bu}{\mbox{\bf u}}
\newcommand{\bv}{\mbox{\bf v}}
\newcommand{\bw}{\mbox{\bf w}}
\newcommand{\bx}{\mbox{\bf x}}
\newcommand{\by}{\mbox{\bf y}}
\newcommand{\bz}{\mbox{\bf z}}
\newcommand{\bA}{\mbox{\bf A}}
\newcommand{\bB}{\mbox{\bf B}}

\newcommand{\bD}{\mbox{\bf D}}

\newcommand{\bG}{\mbox{\bf G}}

\newcommand{\bR}{\mbox{\bf R}}

\newcommand{\bX}{\mbox{\bf X}}

\newcommand{\calS}{{\cal S}}
\newcommand{\bY}{\mbox{\bf Y}}

\newcommand{\bone}{\mbox{\bf 1}}
\newcommand{\bzero}{\mbox{\bf 0}}

\newcommand{\bbeta}{\mbox{\boldmath $\beta$}}
\newcommand{\btheta}{\mbox{\boldmath $\theta$}}
\newcommand{\bgamma}{\mbox{\boldmath $\gamma$}}
\newcommand{\bdelta}{\mbox{\boldmath $\delta$}}
\newcommand{\bPsi}{\mbox{\boldmath $\Psi$}}
\newcommand{\bet}{\mbox{\boldmath $\eta$}}
\newcommand{\bxi}{\mbox{\boldmath $\xi$}}

\newcommand{\bmu}{\mbox{\boldmath $\mu$}}

\newcommand{\hbbeta}{\widehat\bbeta}
\newcommand{\hbtheta}{\widehat\btheta}

\newcommand{\hbeta}{\widehat\beta}

\newcommand{\var}{\mathrm{var}}
\newcommand{\cov}{\mathrm{cov}}

\newcommand{\Sig}{\mathbf{\Sigma}}

\newcommand{\tr}{\mathrm{tr}}
\newcommand{\diag}{\mathrm{diag}}

\newcommand{\sgn}{\mathrm{sgn}}
\newcommand{\supp}{\mathrm{supp}}
\def\etal{\textit{et al.} }
\def\t{^T}
\def\toD{\overset{\mathscr{D}}{\longrightarrow}}

\begin{document}

\title{Non-Concave Penalized Likelihood with NP-Dimensionality} %
\date{September 1, 2009}
\author{Jianqing Fan and Jinchi Lv %
\thanks{Jianqing Fan is Frederick L. Moore '18 Professor of Finance, Department of Operations Research and Financial Engineering, Princeton University, Princeton, NJ 08544, USA (e-mail: jqfan@princeton.edu). %
Jinchi Lv is Assistant Professor of Statistics, Information and Operations Management Department, Marshall School of Business, University of Southern California, Los Angeles, CA 90089, USA (e-mail: jinchilv@marshall.usc.edu). %
Fan's research was partially supported by NSF Grants DMS-0704337 and DMS-0714554 and NIH Grant
R01-GM072611. Lv's research was partially supported by NSF Grant DMS-0806030 and 2008 Zumberge Individual Award from USC's James H. Zumberge Faculty Research and Innovation Fund.
}
\\
Princeton University and University of Southern California\\
}

\maketitle

\begin{abstract}
Penalized likelihood methods are fundamental to ultra-high
dimensional variable selection. How high dimensionality
such methods can handle remains largely unknown. In this paper, we show
that in the context of generalized linear models, such methods
possess model selection consistency with oracle properties even
for dimensionality of Non-Polynomial (NP) order of sample size,
for a class of penalized likelihood approaches using folded-concave penalty functions, which were introduced to ameliorate the bias problems of convex penalty
functions. This fills a long-standing gap in the literature where
the dimensionality is allowed to grow slowly with the sample size.
Our results are also applicable to penalized likelihood with the
$L_1$-penalty, which is a convex function at the boundary of
the class of folded-concave penalty functions under consideration. The coordinate
optimization is implemented for finding the solution paths, whose
performance is evaluated by a few simulation examples and the real
data analysis.
\end{abstract}

\textit{Running title}: Non-Concave Penalized Likelihood


\section{Introduction} \label{Sec1}
The analysis of data sets with the number of variables $p$
comparable to or much larger than the sample size $n$ frequently
arises nowadays in many fields ranging from genomics and health
sciences to economics and machine learning. The data that we collect
is usually of the type $(y_i, x_{i1}, \cdots, x_{ip})_{i = 1}^n$,
where the $y_i$'s are $n$ independent observations of the response
variable $Y$ given its covariates, or explanatory
variables, $(x_{i1}, \cdots, x_{ip})\t$. Generalized linear models
(GLMs) provide a flexible parametric approach to estimating the
covariate effects (McCullagh and Nelder, 1989). In this paper we
consider the variable selection problem of Non-Polynomial (NP)
dimensionality in the context of GLMs.  By NP-dimensionality we
mean that $\log p = O(n^a)$ for some $a \in (0, 1)$. See Fan and Lv (2009) for an overview of recent developments in high dimensional variable selection.

We denote by $\bX = (\bx_1, \cdots, \bx_p)$ the $n \times p$ design
matrix with $\bx_j = (x_{1j}, \cdots, x_{nj})\t$, $j = 1, \cdots, p$
and $\by = (y_1, \cdots, y_n)\t$ the $n$-dimensional response
vector. Throughout the paper we consider deterministic design
matrix. With a canonical link, the conditional distribution of $\by$
given $\bX$ belongs to the canonical exponential family, having
the following density function with respect to some fixed measure
\begin{equation} \label{001}
f_n(\by; \bX, \bbeta) \equiv \prod_{i = 1}^n f_0(y_i; \theta_i) = \prod_{i = 1}^n \left\{c(y_i) \exp\left[\frac{y_i \theta_i - b(\theta_i)}{\aphi}\right]\right\},
\end{equation}
where $\bbeta = (\beta_1, \cdots, \beta_p)\t$ is an unknown
$p$-dimensional vector of regression coefficients, $\{f_0(y; \theta): \theta \in \mathbf{R}\}$ is a family of distributions in
the regular exponential family with dispersion parameter $\aphi \in (0, \infty)$, and
$(\theta_1, \cdots, \theta_n)\t = \bX \bbeta$. As is common in GLM,
the function $b(\theta)$ is implicitly assumed to be twice
continuously differentiable with $b''(\theta)$ always positive. In
the sparse modeling, we assume that majority of the true
regression coefficients $\bbeta_0 = (\beta_{0, 1}, \cdots,
\beta_{0, p})\t$ are exactly zero. Without
loss of generality, assume that $\bbeta_0 = (\bbeta_1\t,
\bbeta_2\t)\t$ with each component of $\bbeta_1$ nonzero and
$\bbeta_2 = \bzero$. Hereafter we refer to the support $\supp(\bbeta_0) = \{1,
\cdots, s\}$ as the true underlying sparse model of the indices. Variable selection
aims at locating those predictors $\bx_j$ with nonzero $\beta_{0,
j}$ and giving an effective estimate of $\bbeta_1$.

In view of (\ref{001}), the log-likelihood $\log f_n(\by; \bX,
\bbeta)$ of the sample is given, up to an affine transformation, by
\begin{equation} \label{002}
\ell_n(\bbeta) = n^{-1}  \left[\by\t \bX \bbeta - \bone\t \bb(\bX
\bbeta)\right],
\end{equation}
where $\bb(\btheta) = (b(\theta_1), \cdots, b(\theta_n))\t$ for $\btheta = (\theta_1, \cdots, \theta_n)\t$. We consider the following penalized likelihood
\begin{equation} \label{004}
Q_n(\bbeta) = \ell_n(\bbeta) - \sum_{j = 1}^p p_{\lambda_n}(|\beta_j|),
\end{equation}
where $p_{\lambda}(\cdot)$ is a penalty function and $\lambda_n \geq 0$ is a regularization parameter.

In a pioneering paper, Fan and Li (2001) build the theoretical
foundation of nonconcave penalized likelihood for variable
selection. The penalty functions that they used are not any
nonconvex functions, but really the folded-concave functions. For
this reason, we will call them more precisely folded-concave
penalties. The paper also introduces the oracle property for model
selection. An estimator $\hbbeta = (\hbbeta_1\t, \hbbeta_2\t)\t$ is
said to have the oracle property (Fan and Li, 2001) if it enjoys the
model selection consistency in the sense of $\hbbeta_2 = \bzero$
with probability tending to 1 as $n \rightarrow \infty$, and it
attains an information bound mimicking that of the oracle estimator,
where $\hbbeta_1$ is a subvector of $\hbbeta$ formed by its first
$s$ components and the oracle knew the true model $\supp(\bbeta_0) =
\{1, \cdots, s\}$ ahead of time. Fan and Li (2001) study the oracle
properties of non-concave penalized likelihood estimators in the
finite-dimensional setting. Their results were extended later by Fan
and Peng (2004) to the setting of $p = o(n^{1/5})$ or $o(n^{1/3})$
in a general likelihood framework. The question of how large $p$ can
be so that similar oracle properties continue to hold arises
naturally. Can the penalized likelihood methods be applicable to
NP-dimensional variable selection problems?  This paper gives an
affirmative answer.

Numerous efforts have lately been devoted to studying the properties of
variable selection with ultra-high dimensionality and significant
progress has been made. Meinshausen and B\"{u}hlmann (2006), Zhao and
Yu (2006), and Zhang and Huang (2008) investigate the issue of model
selection consistency for LASSO under different setups when the number of
variables is of a greater order than the sample size. Candes and Tao
(2007) introduce the Dantzig selector to handle the NP-dimensional
variable selection problem, which was shown to behave similarly to
Lasso by Bickel \etal (2009). Zhang (2009) is among the first to study the non-convex penalized least-squares estimator with NP-dimensionality and demonstrates its advantages over LASSO.  He also develops the PLUS algorithm to find the solution path that has the desired sampling properties.
Fan and Lv (2008) and Huang \etal
(2008) introduce the independence screening procedure to reduce the
dimensionality in the context of least-squares. The former
establishes the sure screening property with NP-dimensionality
and the latter also studies the bridge regression, a folded-concave
penalty approach. Fan and Fan (2008) investigate the impact of dimensionality
on ultra-high dimensional classification and establish an oracle
property for features annealed independence
rules. Lv and Fan (2009) make important connections between model selection and sparse
recovery using folded-concave penalties and establish a nonasymptotic weak oracle property for the penalized least squares estimator with NP-dimensionality. There are also a number of
important papers on establishing the oracle inequalities for
penalized empirical risk minimization. For example, Bunea \etal (2007)
establish sparsity oracle inequalities for the Lasso under quadratic
loss in the context of least-squares; van de Geer (2008) obtains a
nonasymptotic oracle inequality for the empirical risk minimizer
with the $L_1$-penalty in the context of GLMs; Koltchinskii (2008) proves oracle inequalities for penalized least squares with entropy penalization.

The penalization methods are also widely used in covariance matrix estimation.  This has been studied by a number of authors on the estimation of sparse covariance matrix, sparse precision matrix, and sparse Cholesky decomposition, using the Gaussian likelihood or pseudo-likelihood.  See, for example,  Huang \etal (2006), Meinshausen and B\"{u}hlmann (2006), Levina \etal (2008), Rothman \etal (2008), and Lam and Fan (2009), among others.  For these more specific models, stronger results can be obtained.

The rest of the paper is organized as follows. In Section \ref{Sec2}, we
discuss the choice of penalty functions and characterize the
non-concave penalized likelihood estimator and its global optimality. We study the
nonasymptotic weak oracle properties and oracle properties of
non-concave penalized likelihood estimator in Sections \ref{Sec3} and \ref{Sec4},
respectively. Section \ref{Sec5} discusses algorithms for solving
regularization problems with concave penalties including the SCAD. In
Section \ref{Sec6}, we present three numerical examples using both simulated and
real data sets.  We provide some discussions of our results and their implications in Section \ref{Sec7}. Proofs are presented in Section \ref{Sec8}.
Technical details are relegated to the Appendix.

\section{Non-concave penalized likelihood estimation} \label{Sec2}
In this section we discuss the choice of penalty functions in
regularization methods and characterize the non-concave penalized
likelihood estimator as well as its global optimality.

\subsection{Penalty function} \label{Sec2.1}
For any penalty function $p_\lambda(\cdot)$, we let $\rho(t; \lambda) = \lambda^{-1} p_\lambda(t)$. For simplicity, we will drop its dependence on $\lambda$ and write $\rho(t; \lambda)$ as $\rho(t)$ when there is no confusion. Many penalty functions have been proposed in the literature for regularization. For example, the best subset selection amounts to using the $L_0$ penalty. The ridge regression uses the $L_2$ penalty. The $L_q$ penalty $\rho(t) = t^q$ for $q \in (0, 2)$  bridges these two cases (Frank and Friedman, 1993). Breiman (1995) introduces the non-negative garrote for shrinkage estimation and variable selection. Lasso (Tibshirani, 1996) uses the $L_1$-penalized least squares. The SCAD penalty (Fan, 1997; Fan and Li, 2001) is the function whose derivative is given by
\begin{equation} \label{003}
p_\lambda'(t) = \lambda \left\{I\left(t \leq \lambda\right) + \frac{\left(a \lambda - t\right)_+}{\left(a - 1\right) \lambda} I\left(t > \lambda\right)\right\}, \quad t \ge 0, \text{ for some } a > 2,
\end{equation}
where often $a = 3.7$ is used, and MCP (Zhang, 2009) is defined
through $p_\lambda'(t) =\left(a \lambda - t\right)_+/a$. Clearly the
SCAD penalty takes off at the origin as the $L_1$ penalty and then
levels off, and MCP translates the flat part of the derivative of
SCAD to the origin. A family of folded concave penalties that bridge the
$L_0$ and $L_1$ penalties were studied by Lv and Fan (2009).

Hereafter we consider penalty functions $p_\lambda(\cdot)$ that satisfy the following condition:

\begin{assumption} \label{con1}
$\rho(t; \lambda)$ is increasing and concave in $t \in [0, \infty)$, and has
a continuous derivative $\rho'(t; \lambda)$ with $\rho'(0+; \lambda)> 0$.  In addition, $\rho'(t; \lambda)$ is increasing
in $\lambda \in (0, \infty)$ and $\rho'(0+; \lambda)$ is independent of
$\lambda$.
\end{assumption}

The above class of penalty functions has been considered by Lv and Fan (2009). Clearly the $L_1$ penalty is a convex function that falls at the
boundary of the class of penalty functions satisfying Condition
\ref{con1}. Fan and Li (2001) advocate penalty functions that give
estimators with three desired properties: unbiasedness, sparsity and
continuity, and provide insights into them (see also Antoniadis and
Fan, 2001). Both SCAD and MCP with $a \geq 1$ satisfy Condition
\ref{con1} and the above three properties simultaneously. The $L_1$
penalty also satisfies Condition \ref{con1} as well as the sparsity
and continuity, but it does not enjoy the unbiasedness, since its
derivative is identically one on $[0, \infty)$ with the derivative at zero understood as the right derivative. However, our results
are applicable to the $L_1$-penalized regression. Condition \ref{con1} is needed
for establishing the oracle properties of non-concave penalized likelihood estimator.

\subsection{Non-concave penalized likelihood estimator} \label{Sec2.2}
It is generally difficult to study the global maximizer of the penalized likelihood analytically without concavity. As is common in the literature, we study the behavior of local maximizers.

We introduce some notation to simplify our presentation. For any $\btheta = (\theta_1, \cdots, \theta_n)\t \in \mathbf{R}^n$, define
\begin{equation} \label{105}
\bmu(\btheta) = (b'(\theta_1), \cdots, b'(\theta_n))\t \ \text{ and } \
\Sig(\btheta) = \diag\{b''(\theta_1), \cdots, b''(\theta_n)\}.
\end{equation}
It is known that the $n$-dimensional response vector $\by$ following
the distribution in (\ref{001}) has mean vector $\bmu(\btheta)$ and
covariance matrix $\aphi \Sig(\btheta)$, where $\btheta = \bX
\bbeta$. Let $\bar{\rho}(t) = \sgn(t) \rho'(|t|)$, $t \in
\mathbf{R}$ and $\bar{\rho}(\bv) = (\bar{\rho}(v_1), \cdots,
\bar{\rho}(v_q))\t$, $\bv = (v_1, \cdots, v_q)\t$, where $\sgn$
denotes the sign function. We denote by $\|\cdot\|_q$ the $L_q$ norm of a vector or matrix for $q \in [0, \infty]$. Following Zhang (2009), define the local
concavity of the penalty $\rho$ at $\bv = (v_1, \cdots, v_q)\t \in
\mathbf{R}^q$ with $\|\bv\|_0 = q$ as
\begin{equation} \label{016}
\kappa(\rho; \bv) = \lim_{\epsilon \rightarrow 0+} \max_{1 \leq j \leq q} \sup_{t_1 <
 t_2 \in (|v_j| - \epsilon, |v_j| + \epsilon)} -\frac{\rho'(t_2) - \rho'(t_1)}{t_2 - t_1}.
\end{equation}
By the concavity of $\rho$ in Condition \ref{con1}, we have
$\kappa(\rho; \bv) \geq 0$. It is easy to show by the mean-value
theorem that $\kappa(\rho; \bv) = \max_{1 \leq j \leq q}
-\rho''(|v_j|)$ provided that the second derivative of $\rho$ is
continuous. For the SCAD penalty, $\kappa(\rho; \bv) = 0$ unless some component of $|\bv|$ takes values in $[\lambda, a \lambda]$. In the latter case, $\kappa(\rho; \bv) = (a-1)^{-1} \lambda^{-1}$.

Throughout the paper, we use $\lambda_{\min}(\cdot)$ and
$\lambda_{\max}(\cdot)$ to represent the smallest and largest
eigenvalues of a symmetric matrix, respectively.

The following theorem gives a sufficient condition on the strict local maximizer of the penalized likelihood $Q_n(\bbeta)$ in (\ref{004}) (see Lv and Fan (2009) for the case of penalized least squares).

\begin{theorem}[Characterization of PMLE] \label{T1}
Assume that $p_\lambda$ satisfies Condition \ref{con1}. Then $\hbbeta \in \mathbf{R}^p$ is a strict local maximizer of the non-concave penalized likelihood $Q_n(\bbeta)$ defined by (\ref{004}) if
\begin{align}
\label{021} & \bX_1\t \by - \bX_1\t \bmu(\hbtheta) -  n \lambda_n \bar{\rho}(\hbbeta_1) = \bzero, \\
\label{022} & \|\bz\|_\infty < \rho'(0+),\\
\label{023} & \lambda_{\min} \left[\bX_1\t \Sig\left(\hbtheta\right)
\bX_1\right] >  n \lambda_n \kappa(\rho; \hbbeta_1),
\end{align}
where $\bX_1$ and $\bX_2$ respectively denote the submatrices of
$\bX$ formed by columns in $\supp(\hbbeta)$ and its complement,
$\hbtheta = \bX \hbbeta$, $\hbbeta_1$ is a subvector of $\hbbeta$
formed by all nonzero components, and $\bz = ( n \lambda_n)^{-1}
\bX_2\t [\by - \bmu(\hbtheta)]$.  On the other hand, if $\hbbeta$ is
a local maximizer of $Q_n(\bbeta)$, then it must satisfy (\ref{021})
-- (\ref{023}) with strict inequalities replaced by nonstrict
inequalities.
\end{theorem}

There is only a tiny gap (nonstrict versus strict inequalities)
between the necessary condition for local maximizer and sufficient
condition for strict local maximizer.
Conditions (\ref{021}) and (\ref{023}) ensure that $\hbbeta$ is a
strict local maximizer of (\ref{004}) when constrained on the
$\|\hbbeta\|_0$-dimensional subspace $\{\bbeta \in \mathbf{R}^p:
\bbeta_c = \bzero\}$ of $\mathbf{R}^p$, where $\bbeta_c$ denotes the
subvector of $\bbeta$ formed by components in the complement of
$\supp(\hbbeta)$. Condition (\ref{022}) makes sure that the sparse
vector $\hbbeta$ is indeed a strict local maximizer of (\ref{004})
on the whole space $\mathbf{R}^p$.

When $\rho$ is the $L_1$ penalty, the penalized likelihood function $Q_n(\bbeta)$ in (\ref{004}) is concave in $\bbeta$. Then the classical convex optimization theory applies to show that $\hbbeta = (\hbeta_1, \cdots, \hbeta_p)\t$ is a global maximizer if and only if there exists a subgradient $\bz \in \partial L_1(\hbbeta)$ such that
\begin{equation} \label{SGEq1}
\bX\t \by - \bX\t \bmu(\hbtheta) -  n \lambda_n \bz = \bzero,
\end{equation}
that is, it satisfies the Karush-Kuhn-Tucker (KKT) conditions, where
the subdifferential of the $L_1$
penalty is given by $\partial L_1(\hbbeta) = \{ \bz = (z_1, \cdots,
z_p)\t \in \mathbf{R}^p: z_j = \sgn(\hbeta_j) \text{ for } \hbeta_j
\neq 0 \text{ and } z_j \in [-1, 1] \text{ otherwise}\}$. Thus
condition (\ref{SGEq1}) reduces to (\ref{021}) and (\ref{022}) with
strict inequality replaced by nonstrict
inequality.  Since $\kappa(\rho; \bv) = 0$ for the $L_1$-penalty, condition (\ref{023}) holds provided that $\bX_1\t \Sig(\hbtheta) \bX_1$ is nonsingular.
However, to ensure that $\hbbeta$ is the strict maximizer we need
the strict inequality in (\ref{022}).

\subsection{Global optimality} \label{Sec2.3}

It is a natural question of when the non-concave penalized maximum likelihood estimator (NCPMLE) $\hbbeta$ is a global maximizer of the penalized likelihood $Q_n(\bbeta)$. We characterize such a property from two perspectives.

\subsubsection{Global optimality} \label{Sec2.3.1}
Assume that the $n \times p$ design matrix $\bX$ has full column rank $p$. This implies that $p \leq n$. Since $b''(\theta)$ is always positive, it is easy to show that the Hessian matrix of $-\ell_n(\bbeta)$ is always positive definite, which entails that the log-likelihood function $\ell_n(\bbeta)$ is strictly concave in $\bbeta$. Thus there exists a unique maximizer $\bbeta_*$ of $\ell_n(\bbeta)$.
Let $\mathcal{L}_c = \{\bbeta \in \mathbf{R}^p: \ell_n(\bbeta) \geq c \}$ be a sublevel set of $-\ell_n(\bbeta)$ for some $c < \ell_n(\bzero)$ and
\[
    \kappa(p_{\lambda}) = \sup_{t_1 < t_2 \in (0, \infty)} -\frac{p_{\lambda}'(t_2) - p_{\lambda}'(t_1)}{t_2 - t_1}
\]
be the maximum concavity of the penalty function $p_\lambda$. For the $L_1$ penalty, SCAD, and MCP, we have $\kappa(p_\lambda) = 0$, $(a - 1)^{-1}$, and
$a^{-1}$, respectively. The following proposition gives a sufficient condition on the global optimality of NCPMLE.

\begin{proposition}[Global optimality] \label{P2}
Assume that $\bX$ has rank $p$ and satisfies
\begin{equation}\label{141}
\min_{\bbeta \in \mathcal{L}_c} \lambda_{\min} \left[n^{-1} \bX\t
\Sig\left(\bX \bbeta\right) \bX\right] \geq  \kappa(p_{\lambda_n}).
\end{equation}
Then the NCPMLE $\hbbeta$ is a global maximizer of the penalized likelihood $Q_n(\bbeta)$ if $\hbbeta \in \mathcal{L}_c$.
\end{proposition}

Note that for penalized least-squares,  (\ref{141}) reduces to
\begin{equation}\label{14a}
  \lambda_{\min} \left(n^{-1} \bX\t  \bX\right) \geq
  \kappa(p_{\lambda_n}).
\end{equation}
This condition holds for sufficiently large $a$ in SCAD and MCP, when the correlation between covariates is not too strong. The latter holds for design matrices constructed by using spline bases to approximate a nonparametric function.  According to Proposition~\ref{P2}, under (\ref{14a}), the penalized least-squares with folded-concave penalty is a global minimum.

The proposition below gives a condition under which the penalty term in (\ref{004}) does not change the global maximizer.  It will be used to derive the condition under which the PMLE is the same as the oracle estimator in Proposition \ref{P4}(b).  Here for simplicity we consider the SCAD penalty $p_\lambda$ given by (\ref{003}), and the technical arguments are applicable to other folded-concave penalties as well.

\begin{proposition}[Robustness] \label{P3}
Assume that $\bX$ has rank $p$ with $p = s$ and there exists some $c
< \ell_n(\bzero)$ such that $\min_{\bbeta \in \mathcal{L}_c}
\lambda_{\min} [n^{-1} \bX\t \Sig(\bX \bbeta) \bX] \geq c_0 $ for
some $c_0 > 0$. Then the SCAD penalized likelihood estimator
$\hbbeta$ is the global maximizer and equals $\bbeta_*$ if $\hbbeta \in \mathcal{L}_c$ and $\min_{j = 1}^p |\hbeta_j| > (a + \frac{1}{2 c_0}) \lambda_n$, where $\hbbeta = (\hbeta_1, \cdots, \hbeta_p)\t$.
\end{proposition}

\subsubsection{Restricted global optimality} \label{Sec2.3.2}

When $p > n$, it is hard to show the global optimality of a local maximizer. However, we can study the global optimality of the NCPMLE $\hbbeta$ on the union of coordinate subspaces. A subspace of $\mathbf{R}^p$ is called coordinate subspace if it is spanned by a subset of the natural basis $\{\be_1, \cdots, \be_p\}$, where each $\be_j$ is the $p$-vector with $j$-th component 1 and 0 elsewhere. Here each $\be_j$ corresponds to the $j$-th predictor $\bx_j$.  We will investigate the global optimality of $\hbbeta$ on the union $\mathbb{S}_s$ of all $s$-dimensional coordinate subspaces of $\mathbf{R}^p$ in Proposition~\ref{P4}(a).

Of particularly interest is to derive the conditions under which the PMLE is also an oracle estimator, in addition to possessing the above restricted global optimal estimator on $\mathbb{S}_s$.  To this end, we introduce an identifiability condition on the true model $\supp(\bbeta_0)$. The true model is called $\delta$-identifiable for some $\delta > 0$ if
\begin{equation} \label{143}
\max_{\bbeta \in \mathcal{A}_0} \ell_n(\bbeta) - \sup_{\bbeta \in \mathbb{S}_s \setminus \mathcal{A}_0} \ell_n(\bbeta) \geq \delta,
\end{equation}
where $\mathcal{A}_0 = \{(\beta_1, \cdots, \beta_p)\t \in \mathbf{R}^p: \beta_j = 0 \text{ for } j \notin \supp(\bbeta_0)\}$.  In other words, $\supp(\bbeta_0)$ is the best subset of size $s$, with a margin at least $\delta$.  The following proposition is an easy consequence of Propositions \ref{P2} and \ref{P3}.

\begin{proposition}[Global optimality on $\mathbb{S}_s$] \hspace{0.05 in}
\label{P4}
\begin{itemize}
\item[a)]
If the conditions of Proposition \ref{P2} are satisfied for each $n \times (2 s)$ submatrix of $\bX$, then the NCPMLE $\hbbeta$ is a global maximizer of $Q_n(\bbeta)$ on $\mathbb{S}_s$.

\item[b)] Assume that the conditions of Proposition \ref{P3} are satisfied for the $n \times s$ submatrix of $\bX$ formed by columns in $\supp(\bbeta_0)$, the true model is $\delta$-identifiable for some $\delta > \frac{(a + 1) s \lambda_n^2}{2}$, and $\supp(\hbbeta) = \supp(\bbeta_0)$. Then the SCAD penalized likelihood estimator $\hbbeta$ is the global maximizer on $\mathbb{S}_s$ and equals to the oracle maximum likelihood estimator $\bbeta_*$.
\end{itemize}
\end{proposition}

On the event that the PMLE estimator is the same as the oracle estimator, it possesses of course the oracle property.

\section{Nonasymptotic weak oracle properties} \label{Sec3}
In this section we study a nonasymptotic property of the
non-concave penalized likelihood estimator $\hbbeta$, called the
weak oracle property introduced by Lv and Fan (2009) in the setting
of penalized least squares. The weak oracle property means
sparsity in the sense of $\hbbeta_2 = \bzero$ with probability
tending to 1 as $n \rightarrow \infty$, and consistency under the
$L_{\infty}$ loss, where $\hbbeta =(\hbbeta_1\t, \hbbeta_2\t)\t$ and
$\hbbeta_1$ is a subvector of $\hbbeta$ formed by components in
$\supp(\bbeta_0) = \{1, \cdots, s\}$. This property is weaker than
the oracle property introduced by Fan and Li (2001).

\subsection{Regularity conditions} \label{Sec3.1}
As mentioned before, we condition on the design matrix $\bX$ and use the $p_\lambda$ penalty in the class satisfying Condition \ref{con1}. Let $\bX_1$ and $\bX_2$ respectively be the submatrices of the $n \times p$ design matrix $\bX = (\bx_1, \cdots, \bx_p)$ formed by columns in $\supp(\bbeta_0)$ and its complement, and
$\btheta_0 = \bX \bbeta_0$. To simplify the presentation, we assume without loss of generality
that each covariate $\bx_j$ has been standardized so that $\|\bx_j\|_2 =
\sqrt{n}$. If the covariates have not been standardized, the results still hold with $\|\bx_j\|_2$ assumed to be in the order of $\sqrt{n}$. Let
\begin{equation} \label{138}
d_n = 2^{-1} \min\left\{\left|\beta_{0, j}\right|: \beta_{0, j} \neq 0\right\}
\end{equation}
be half of the minimum signal. We make the following assumptions on the
design matrix and the distribution of the response.

Let $\{b_s\}$ be a diverging sequence of positive numbers that depends on the nonsparsity size $s$ and hence depends on $n$. Recall that $\bbeta_1$ is the non-vanishing components of the true parameter $\bbeta_0$.
\begin{assumption} \label{con2}
The design matrix $\bX$ satisfies
\begin{align} \label{007}
& \left\|\left[\bX_1\t \Sig\left(\btheta_0\right) \bX_1\right]^{-1}\right\|_\infty
= O(b_s n^{-1}), \\
\label{008}
& \left\|\bX_2\t \Sig\left(\btheta_0\right) \bX_1 \left[\bX_1\t \Sig\left(\btheta_0\right)
\bX_1\right]^{-1}\right\|_\infty \leq \min\left\{C \frac{\rho'(0+)}{\rho'(d_n)}, O(n^{\alpha_1})\right\}, \\
\label{009}
& \max_{\bdelta \in \mathcal{N}_0} \max\nolimits_{j = 1}^p \lambda_{\max}\left[\bX_1\t \diag\left\{\left|\bx_j\right| \circ \left|\bmu''\left(\bX_1 \bdelta\right)\right|\right\} \bX_1\right] = O(n),
\end{align}
where the $L_\infty$ norm of a matrix is the maximum of the $L_1$ norm of each row, $C \in (0, 1)$, $\alpha_1 \in [0, 1/2]$, $\mathcal{N}_0 = \{\bdelta \in \mathbf{R}^s: \|\bdelta -
\bbeta_1\|_\infty \leq d_n\}$, the derivative is taken componentwise, and $\circ$ denotes the Hadamard (componentwise) product.
\end{assumption}

Here and below, $\rho$ is associated with regularization parameter
$\lambda_n$ satisfying (\ref{104}) unless specified
otherwise. For the classical Gaussian linear regression model, we
have $\bmu(\btheta) = \btheta$ and $\Sig(\btheta) =  I_n$. In
this case, since we will assume that $s \ll n$, condition (\ref{007}) usually holds with $b_s = 1$ if the covariates are nearly uncorrelated.  In fact, Wainwright (2009) shows that $\|[\bX_1^T\bX_1]^{-1}\|_\infty = O_P(n^{-1})$ if the rows of $\bX_1$ are i.i.d. Gaussian vectors with $\| [E \bX_1^T\bX_1]^{-1}\|_\infty = O_P(n^{-1})$.  In general, since
$$
\|[\bX_1^T\bX_1]^{-1}\|_\infty \leq \sqrt{s}/\lambda_{\min} (\bX_1^T\bX_1),
$$
we can take $b_s = s^{1/2}$ if $\lambda_{\min} (\bX_1^T\bX_1)^{-1} = O(n^{-1})$.  More generally, (\ref{007}) can be bounded as
$$
\left\|\left[\bX_1\t \Sig\left(\btheta_0\right) \bX_1\right]^{-1}\right\|_\infty = d^{-1} \left\|\left[\bX_{1,S}\t \bX_{1,S}\right]^{-1}\right\|_\infty
$$
and the above remark for the multiple regression model applies to the submatrix $\bX_{1, S}$, which consists of rows of the samples with $b''(\theta_i) > d$ for some $d>0$.

The left hand side of (\ref{008}) is the multiple regression coefficients
of each unimportant variable in $\bX_2$ on $\bX_1$, using the weighted least squares with weights $\{b''(\theta_i)\}$. Condition (\ref{008}) controls the uniform growth rate of the $L_1$-norm of these multiple regression coefficients, a notion of weak correlation between $\bX_1$ and $\bX_2$. If each element of the multiple regression coefficients is of order $O(1)$, then the $L_1$ norm is of order $O(s)$.  Hence, we can handle the non-sparse dimensionality $s = O(n^{\alpha_1})$, by (\ref{008}), as long as the first term in (\ref{008}) dominates, which occurs for SCAD type of penalty with $d_n \gg \lambda_n$.  Of course, the actual dimensionality can be higher or lower, depending on the correlation between $\bX_1$ and $\bX_2$, but for finite non-sparse dimensionality $s=O(1)$, (\ref{008}) is usually satisfied.
When a folded-concave penalty is used,
the upper bound on the right hand side of (\ref{008}) can grow to
$\infty$ at rate $O(n^{\alpha_1})$.  In  contrast,
when the $L_1$ penalty is used, the upper bound in (\ref{008}) is more restrictive, requiring uniformly less than 1.  This condition is the same as the strong irrepresentable condition of Zhao and Yu (2006) for the consistency of the LASSO estimator, namely $\|\bX_2\t
\bX_1 (\bX_1\t \bX_1)^{-1}\|_\infty \leq C$.  It is a
drawback of the $L_1$ penalty.

For the Gaussian linear regression
model, condition (\ref{009}) holds automatically.

We now choose the regularization parameter $\lambda_n$ and introduce Condition \ref{con3}. We will assume that half of the minimum signal $d_n \geq n^{-\gamma} \log n$ for some $\gamma \in (0, 1/2]$. Take $\lambda_n$ satisfying
\begin{eqnarray} \label{104}
p_{\lambda_n}'(d_n) = o(b_s^{-1} n^{-\gamma} \log n) \quad \text{and} \quad \lambda_n \gg n^{-\alpha} (\log n)^{2},
\end{eqnarray}
where $\alpha = \min(\frac{1}{2}, 2 \gamma - \alpha_0) - \alpha_1$ and $b_s$ is associated with the nonsparsity size $s = O(n^{\alpha_0})$.

\begin{assumption} \label{con3}
Assume that $d_n \geq n^{-\gamma} \log n$ and
$b_s = o\{\min(n^{1/2 -\gamma} \sqrt{\log n}, \\ s^{-1} n^\gamma / \log n)\}$.  In addition, assume that $\lambda_n$ satisfies (\ref{104}) and ${\lambda}_n \kappa_0 = o(\tau_0)$, where $\kappa_0 = \max_{\bdelta \in \mathcal{N}_0} \kappa(\rho; \bdelta)$ and $\tau_0 = \min_{\bdelta \in \mathcal{N}_0} \lambda_{\min}[n^{-1} \bX_1\t \Sig(\bX_1 \bdelta) \bX_1]$, and that $\max_{j = 1}^p \|\bx_j\|_\infty = o(n^\alpha/\sqrt{\log n})$ if the responses are unbounded.
\end{assumption}

The condition that ${\lambda}_n \kappa_0 = o(\tau_0)$, is needed to ensure condition (\ref{023}).  The condition always holds when $\kappa_0 = 0$ and is satisfied for the SCAD type of penalty when $d_n \gg \lambda_n$.

In view of (\ref{021}) and (\ref{022}), to study the non-concave penalized likelihood
estimator $\hbbeta$ we need to analyze the deviation of the
$p$-dimensional random vector $\bX\t \bY$ from its mean $\bX\t
\bmu(\btheta_0)$, where $\bY = (Y_1, \cdots, Y_n)\t$ denotes the $n$-dimensional random
response vector in the GLM (\ref{001}). The following proposition, whose proof is given in Section~\ref{Sec8.5}, characterizes such
deviation for the case of bounded responses and the case of unbounded
responses satisfying a moment condition, respectively.

\begin{proposition}[Deviation] \label{P1}
Let $\bY = (Y_1, \cdots, Y_n)\t$ be the $n$-dimensional independent random response vector and $\ba \in \mathbf{R}^n$. Then
\begin{itemize}
\item[a)] If $Y_1, \cdots, Y_n$ are bounded in $[c, d]$ for some $c, d \in \mathbf{R}$, then for any $\varepsilon \in (0, \infty)$,
    \begin{equation} \label{029}
    P\left(\left|\ba\t \bY - \ba\t \bmu\left(\btheta_0\right)\right| > \varepsilon\right) \leq 2 \exp\left[-\frac{2 \varepsilon^2}{\left\|\ba\right\|_2^2 (d - c)^2}\right].
    \end{equation}

\item[b)] If $Y_1, \cdots, Y_n$ are unbounded and there exist some $M, v_0 \in (0, \infty)$ such that
    \begin{equation} \label{010}
    \max_{i = 1, \cdots, n} E \left\{\exp\left[\frac{\left|Y_i - b'\left(\theta_{0, i}\right)\right|}{M}\right] - 1 - \frac{\left|Y_i - b'\left(\theta_{0, i}\right)\right|}{M}\right\} M^2 \leq \frac{v_0}{2}
    \end{equation}
    with $(\theta_{0, 1}, \cdots, \theta_{0, n})\t = \btheta_0$, then for any $\varepsilon \in (0, \infty)$,
    \begin{equation} \label{030}
    P\left(\left|\ba\t \bY - \ba\t \bmu\left(\btheta_0\right)\right| > \varepsilon\right) \leq 2 \exp\left[-\frac{1}{2} \frac{\varepsilon^2}{\left\|\ba\right\|_2^2 v_0 + \left\|\ba\right\|_\infty M \varepsilon}\right].
    \end{equation}
\end{itemize}
\end{proposition}

In light of (\ref{001}), it is known that for the exponential family, the moment-generating function of $Y_i$ is given by
\[
E \exp\left\{t \left[Y_i - b'\left(\theta_{0, i}\right)\right]\right\} = \exp\left\{\aphi^{-1} \left[b\left(\theta_{0, i} + t \aphi\right) - b\left(\theta_{0, i}\right) - b'\left(\theta_{0, i}\right) t \aphi\right]\right\},
\]
where $\theta_{0, i} + t \aphi$ is in the domain of $b(\cdot)$. Thus the moment condition (\ref{010}) is reasonable. It is easy to show that condition (\ref{010}) holds for the Gaussian linear regression model and for the Poisson regression model with bounded mean responses. Similar probability bounds also hold for sub-Gaussian errors.

We now express the results in Proposition \ref{P1} in a unified
form. For the case of bounded responses, we define $\varphi(\varepsilon) = 2
e^{-c_1 \varepsilon^2}$ for $\varepsilon \in (0, \infty)$, where $c_1 = 2/(d - c)^2$.
For the case of unbounded
responses satisfying the moment condition (\ref{010}), we define
$\varphi(\varepsilon) = 2 e^{-c_1 \varepsilon^2}$, where $c_1 = 1/(2 v_0 + 2 M)$. Then the exponential bounds in
(\ref{029}) and (\ref{030}) can be expressed as
\begin{equation} \label{137}
P\left(\left|\ba\t \bY - \ba\t \bmu\left(\btheta_0\right)\right| >
  \left\|\ba\right\|_2 \varepsilon \right ) \leq \varphi(\varepsilon),
\end{equation}
where $\varepsilon \in (0, \infty)$ if the responses are bounded and $\varepsilon \in (0, \|\ba\|_2 /\|\ba\|_\infty]$ if the responses are unbounded.

\subsection{Weak oracle properties} \label{Sec3.2}

\begin{theorem}[Weak oracle property] \label{T2}
Assume that Conditions \ref{con1}--\ref{con3} and the probability
bound (\ref{137}) are satisfied, $s = o(n)$, and $\log p = O(n^{1 - 2
\alpha})$. Then there exists a non-concave penalized likelihood
estimator $\hbbeta$ such that for sufficiently large $n$, with
probability at least $1 - 2 [s n^{-1} + (p - s) e^{-n^{1 - 2 \alpha}
\log n}]$, $\hbbeta = (\hbbeta_1\t, \hbbeta_2\t)\t$ satisfies:
\begin{itemize}
\item[a)] \emph{(Sparsity)}. $\hbbeta_2 = \bzero$;

\item[b)] \emph{($L_\infty$ loss)}. $\|\hbbeta_1 - \bbeta_1\|_\infty = O(n^{-\gamma} \log n)$,
\end{itemize}
where $\hbbeta_1$ and $\bbeta_1$ are respectively the subvectors of
$\hbbeta$ and $\bbeta_0$ formed by components in $\supp(\bbeta_0)$.
\end{theorem}

Under the given regularity conditions, the dimensionality $p$ is
allowed to grow up to exponentially fast with the sample size $n$.
The growth rate of $\log p$ is controlled by $1 - 2 \alpha$. It also
enters the nonasymptotic probability bound. This probability tends
to 1 under our technical assumptions.  From the proof of Theorem \ref{T2}, we see that
with asymptotic probability one, the $L_\infty$ estimation loss of
the non-concave penalized likelihood estimator $\hbbeta$ is bounded
from above by three terms (see (\ref{140})), where the second
term $b_s \lambda_n \rho'(d_n)/\rho'(0+)$ is associated with
the penalty function $\rho$. For the $L_1$ penalty, the ratio
$\rho'(d_n)/\rho'(0+)$ is equal to one, and for other concave
penalties, it can be (much) smaller than one. This is in line with
the fact shown by Fan and Li (2001) that concave penalties can
reduce the biases of estimates. Under the specific setting of penalized least squares, the above weak oracle property is slightly different from that of Lv and Fan (2009).

The value of $\gamma$ can be taken as large as $1/2$ for concave penalties.  In this case, the dimensionality that the penalized least-squares can handle is as high as $\log p = O(n^{2\alpha_1})$ when $\alpha_0 \leq 1/2$, which is usually smaller than that for the case of $\gamma < \frac{1}{4} + \frac{\alpha_0}{2}$. The large value of $\gamma$ puts more stringent condition on the design matrix.   To see this, Condition~\ref{con3} entails that $b_s = o(\sqrt{\log n})$ and hence (\ref{007}) becomes tighter.

In the classical setting of $\gamma = 1/2$, the consistency rate of $\hbbeta$ under the $L_2$ norm becomes $O_P(\sqrt{s} n^{-1/2} \log n)$, which is slightly slower than
$O_P(\sqrt{s} n^{-1/2})$. This is because it is derived by using the
$L_\infty$ loss of $\hbbeta$ in Theorem \ref{T2}b). The use of the
$L_\infty$ norm is due to the technical difficulty of proving the
existence of a solution to the nonlinear equation (\ref{021}).

\subsection{Sampling properties of $L_1$-based PMLE} \label{Sec3.3}
When the $L_1$-penalty is applied, the penalized likelihood $Q_n(\bbeta)$ in (\ref{004}) is concave. The local maximizer in Theorems~\ref{T1} and \ref{T2} becomes the global maximizer.   Due to its popularity, we now examine the implications of Theorem~\ref{T2} in the context of penalized least-squares and penalized likelihood.

For the penalized least-squares, Condition~\ref{con2} becomes
\begin{align} \label{007a}
& \left\|\left(\bX_1\t \bX_1\right)^{-1}\right\|_\infty
= O(b_s n^{-1}), \\
\label{008a}
& \left\|\bX_2\t  \bX_1 \left(\bX_1\t
\bX_1\right)^{-1}\right\|_\infty \leq C <1.
\end{align}
Condition (\ref{009}) holds automatically and  Condition (\ref{104}) becomes
\begin{equation}\label{104a}
   \lambda_n = o(b_s^{-1} n^{-\gamma} \log n) \quad \mbox{and} \quad
   \lambda_n \gg n^{-\alpha} (\log n)^{2}.
\end{equation}

As a corollary of Theorem~\ref{T2}, we have

\begin{corollary}[Penalized $L_1$ estimator] \label{C1}
Under Conditions \ref{con2} and \ref{con3} and probability
bound (\ref{137}), if $s = o(n)$ and $\log p = O(n^{1 - 2
\alpha})$, then the penalized $L_1$ likelihood estimator $\hbbeta$ has model selection consistency with rate $\|\hbbeta_1 - \bbeta_1\|_\infty = O(n^{-\gamma} \log n)$.
\end{corollary}

For the penalized least-squares, Corollary~\ref{C1} continues to hold
without normality assumption, as long as probability
bound (\ref{137}) holds.  In this case, the result is stronger than that of Zhao and Yu (2006) and Lv and Fan (2009).

\section{Oracle properties} \label{Sec4}
In this section we study the oracle property (Fan and Li, 2001) of
the non-concave penalized likelihood estimator $\hbbeta$.  We assume
that the nonsparsity size $s \ll n $ and the dimensionality satisfies
$\log p = O(n^\alpha)$ for some $\alpha \in (0, 1/2)$, which is
related to the notation in Section~\ref{Sec3}.   We impose the following
regularity conditions.

\begin{assumption} \label{con4}
The design matrix $\bX$ satisfies
\begin{align} \label{129}
& \min_{\bdelta \in \mathcal{N}_0} \lambda_{\min}\left[\bX_1\t \Sig\left(\bX_1 \bdelta\right) \bX_1\right] \geq c n, \quad \tr[\bX_1\t \Sig(\btheta_0) \bX_1] = O(s n), \\
\label{130}
& \left\|\bX_2\t \Sig\left(\btheta_0\right) \bX_1\right\|_{2, \infty} = O(n),\\
\label{131}
& \max_{\bdelta \in \mathcal{N}_0} \max\nolimits_{j = 1}^p \lambda_{\max}\left[\bX_1\t \diag\left\{\left|\bx_j\right| \circ \left|\bmu''\left(\bX_1 \bdelta\right)\right|\right\} \bX_1\right] = O(n),
\end{align}
where $\mathcal{N}_0 = \{\bdelta \in \mathbf{R}^s: \|\bdelta -
\bbeta_1\|_\infty \leq d_n\}$, $c$ is some positive constant, and $\|\bB\|_{2, \infty} = \max_{\|\bv\|_2 = 1} \|\bB \bv\|_\infty$.
\end{assumption}

\begin{assumption} \label{con5}
Assume that $d_n \gg \lambda_n \gg \max\{(s/n)^{1/2}, n^{(\alpha-1)/2} (\log n)^{1/2}\}$, $p_{\lambda_n}'(d_n) = O(n^{-1/2})$, and $\lambda_n \kappa_0 = o(1)$, where $\kappa_0 = \max_{\bdelta \in
\mathcal{N}_0} \kappa(\rho; \bdelta)$, and in addition that $\max_{j = 1}^p \|\bx_j\|_\infty = o(n^{\frac{1 - \alpha}{2}}/\sqrt{\log n})$ if
the responses are unbounded.
\end{assumption}

Condition~\ref{con4} is generally stronger than Condition~\ref{con2}.  In fact, by $d_n \gg \lambda_n$ in Condition~\ref{con5}, the first condition in (\ref{008}) holds automatically for SCAD type of penalties, since $p_{\lambda_n}'(d_n) = 0$ when $n$ is large enough. Thus Condition \ref{con5} is less restrictive for SCAD-like
penalties, since $\kappa_0 = 0$ for sufficiently large $n$.

However, for the $L_1$ penalty, $\lambda_n =
p_{\lambda_n}'(d_n) = O(n^{-1/2})$ is incompatible with $\lambda_n
\gg (s/n)^{1/2}$. This suggests that the $L_1$
penalized likelihood estimator generally cannot achieve the
consistency rate of $O_P(\sqrt{s} n^{-1/2})$ established in Theorem
\ref{T3} and does not have the oracle property established in
Theorem \ref{T4}, when the dimensionality $p$ is diverging with the
sample size $n$.  In fact, this problem was observed by Fan and Li
(2001) and proved by Zou (2006) even for finite $p$. It still persists with growing dimensionality.

We now state the existence of the NCPMLE and its rate of convergence.  It improves the rate results given by Theorem~\ref{T2}.

\begin{theorem}[Existence of non-concave penalized likelihood estimator] \label{T3}
Assume that Conditions \ref{con1}, \ref{con4} and \ref{con5} and the probability bound (\ref{137}) hold. Then there exists a strict local maximizer $\hbbeta = (\hbbeta_1\t, \hbbeta_2\t)\t$ of the penalized likelihood $Q_n(\bbeta)$ such that $\hbbeta_2 = \bzero$ with probability tending to 1 as $n \rightarrow \infty$ and $\|\hbbeta - \bbeta_0\|_2 = O_P(\sqrt{s} n^{-1/2})$, where $\hbbeta_1$ is a subvector of $\hbbeta$ formed by components in $\supp(\bbeta_0)$.
\end{theorem}

Theorem~\ref{T3} can be thought of as answering the question that given the dimensionality, how strong the minimum signal $d_n$ should be in order for the penalized likelihood estimator to have some nice properties, through Conditions~\ref{con4} and \ref{con5}.  On the other hand, Theorem~\ref{T2} can be thought of as answering the question that
given the strength of the minimum signal $d_n$, how high dimensionality the penalized likelihood methods can handle, through Conditions~\ref{con2} and \ref{con3}. While the details are different, these conditions are related.

To establish the asymptotic normality, we need additional condition, which is related to the Lyapunov condition.

\begin{assumption} \label{con6}
Assume that $p_{\lambda_n}'(d_n) = o(s^{-1/2} n^{-1/2})$, $\max_{i = 1}^n E|Y_i - b'(\theta_{0, i})|^3 = O(1)$, and $\sum_{i = 1}^n (\bz_i\t \bB_n^{-1} \bz_i)^{3/2} \rightarrow 0$ as $n \rightarrow \infty$, where $(Y_1, \cdots, Y_n)\t$ denotes the $n$-dimensional random response vector, $(\theta_{0, 1}, \cdots, \theta_{0, n})\t = \btheta_0$, $\bB_n = \bX_1\t \Sig(\btheta_0) \bX_1$, and $\bX_1 = (\bz_1, \cdots, \bz_n)\t$.
\end{assumption}

\begin{theorem}[Oracle property] \label{T4}
Under the conditions of Theorem~\ref{T3}, if Condition~\ref{con6} holds and $s = o(n^{1/3})$,
then with probability tending to 1 as $n \rightarrow \infty$, the
non-concave penalized likelihood estimator $\hbbeta = (\hbbeta_1\t, \hbbeta_2\t)\t$ in Theorem
\ref{T3} must satisfy:
\begin{itemize}
\item[a)] \emph{(Sparsity)}. $\hbbeta_2 = \bzero$;

\item[b)] \emph{(Asymptotic normality)}.
\[ \bA_n \left[\bX_1\t \Sig\left(\btheta_0\right) \bX_1\right]^{1/2} \left(\hbbeta_1 - \bbeta_1\right) \toD N(\bzero, \aphi \bG), \]
\end{itemize}
where $\bA_n$ is a $q \times s$ matrix such that $\bA_n \bA_n\t \rightarrow \bG$, $\bG$ is a $q \times q$ symmetric positive definite matrix, and $\hbbeta_1$ is a subvector of $\hbbeta$ formed by components in $\supp(\bbeta_0)$.
\end{theorem}

From the proof of Theorem \ref{T4}, we see that for the Gaussian linear regression model, the additional restriction of $s = o(n^{1/3})$ can be relaxed, since the term in (\ref{131}) vanishes in this case.

\section{Implementation} \label{Sec5}
In this section, we discuss algorithms for maximizing the penalized
likelihood $Q_n(\bbeta)$ in (\ref{004}) with concave penalties
including the SCAD. Efficient algorithms for maximizing non-concave
penalized likelihood include the LQA proposed by Fan and Li (2001)
and LLA introduced by Zou and Li (2008). The coordinate optimization
algorithm was used by Fu (1998) and Daubechies \etal (2004) for
penalized least-squares with $L_q$-penalty. This algorithm can also be applied to optimize the group Lasso (Antoniadis and Fan, 2001; Yuan
and Lin, 2006) as shown in Meier \etal (2008) and the penalized
precision matrix estimation in Friedman \etal (2007).

In this paper we employ a path-following algorithm, called the
iterative coordinate ascent (ICA) algorithm. Coordinate optimization
type algorithms are especially appealing for large scale problems
with both $n$ and $p$ large. It successively maximizes $Q_n(\bbeta)$
for regularization parameter $\lambda$ in a decreasing order. ICA
uses the Gauss-Seidel method, i.e., maximizing one coordinate at a
time with successive displacements. Specifically, for each
coordinate within each iteration, ICA uses the second order
approximation of $\ell_n(\bbeta)$ at the $p$-vector from the
previous step along that coordinate and maximizes the univariate
penalized quadratic approximation. It updates each coordinate if the
maximizer of the corresponding univariate penalized quadratic
approximation makes $Q_n(\bbeta)$ strictly increase. Therefore, ICA
algorithm enjoys the ascent property, i.e., the resulting sequence
of $Q_n$ values is increasing for a fixed $\lambda$.

When $\ell_n(\bbeta)$ is quadratic in $\bbeta$, e.g., for the Gaussian linear regression model, the second order approximation in ICA is exact at each step. For any $\bdelta \in \mathbf{R}^p$ and $j \in \{1, \cdots, p\}$, we denote by $\widetilde{\ell}_n(\bbeta; \bdelta, j)$ the second order approximation of $\ell_n(\bbeta)$ at $\bdelta$ along the $j$-th component, and
\begin{equation} \label{144}
  \widetilde{Q}_n(\beta_j; \bdelta, j) = \widetilde{\ell}_n(\bbeta; \bdelta, j) - \sum_{j = 1}^p p_\lambda(|\beta_j|),
\end{equation}
where the subvector of $\bbeta$ with components in $\{1, \cdots, p\}
\setminus \{j\}$ is identical to that of $\bdelta$. Clearly
maximizing $\widetilde{Q}_n(\cdot; \bdelta, j)$ is a univariate
penalized least squares problem, which admits analytical solution for many commonly used penalty functions. See the Appendix for formulae for three popular GLMs.

Pick $\lambda_{\max} \in (0, \infty)$ sufficiently large such that
the maximizer of $Q_n(\bbeta)$ with $\lambda = \lambda_{\max}$ is
$\bzero$,  a decreasing sequence  of regularization parameters
$\{\lambda_1, \cdots, \lambda_K\}$ with $\lambda_1 =
\lambda_{\max}$, and the number of iterations $L$.

\smallskip

\textsc{ICA algorithm}.
\begin{itemize}
\item[1.] Set $k = 1$ and initialize $\hbbeta^{\lambda_0} = \bzero$.

\item[2.] Initialize $\hbbeta^{\lambda_k} = \hbbeta^{\lambda_{k - 1}}$, and set $S = \{1, \cdots, p\}$ and $\ell = 1$.

\item[3.] Successively for $j \in S$, let $\hbeta_j$ be the maximizer of
$\widetilde{Q}_n(\beta_j; \hbbeta^{\lambda_k}, j)$, and update the
$j$-th component of $\hbbeta^{\lambda_k}$ as $\hbeta_j$ if the
updated $\hbbeta^{\lambda_k}$ strictly increases $Q_n(\bbeta)$. Set
$S \leftarrow \supp(\hbbeta^{\lambda_k}) \cup \{j: |z_j| >
\rho'(0+)\}$ and $\ell \leftarrow \ell + 1$, where $(z_1, \cdots,
z_p)\t = ( n \lambda_k)^{-1} \bX\t [\by - \bmu(\bX
\hbbeta^{\lambda_k})]$.

\item[4.] Repeat Step 3 until convergence or $\ell = L + 1$. Set $k \leftarrow k + 1$.

\item[5.] Repeat Steps 2--4 until $k = K + 1$. Return $p$-vectors $\hbbeta^{\lambda_1}, \cdots, \hbbeta^{\lambda_K}$.
\end{itemize}

When we decrease the regularization parameter from $\lambda_k$ to
$\lambda_{k + 1}$, using $\hbbeta^{\lambda_k}$ as an initial value
for $\hbbeta^{\lambda_{k + 1}}$ can speed up the convergence. The
set $S$ is introduced in Step 3 to reduce the computational cost. It
is optional to add $\{j: |z_j| > \rho'(0+)\}$ to the set $S$ in this
step. In practice, we can set a small tolerance level for
convergence. We can also set a level of sparsity for early stopping
if desired models are only those with size up to a certain level.
When the $L_1$ penalty is used, it is known that the choice of
$\lambda =  n^{-1} \|\bX\t [\by - \bmu(\bzero)]\|_\infty$ ensures
that $\bzero$ is the global maximizer of (\ref{004}). In practice,
we can use this value as a proxy for $\lambda_{\max}$. We give the
formulas for three commonly used GLMs and the univariate SCAD
penalized least squares solution in Sections \ref{SecA.1} and \ref{SecA.2}
in the Appendix, respectively.

\section{Numerical examples} \label{Sec6}

\begin{figure} \centering
\begin{center}%
\begin{tabular}
[l]{l}%
{\hspace{-0.45in}\includegraphics[scale=0.75]%
{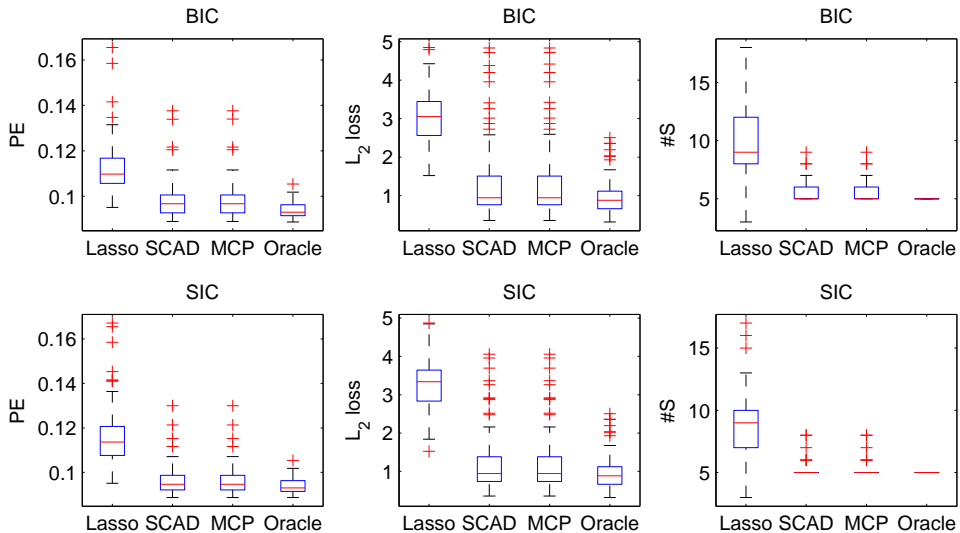}%
}%
\end{tabular}%
\vspace{-0.25in}
\caption{Boxplots of PE, $L_2$ loss, and \#S over $100$ simulations for all methods in logistic regression, where $p = 25$. The $x$-axis represents different methods. Top panel is for BIC and bottom panel is for SIC.}
\label{Fig1}%
\end{center}%
\end{figure}%

\begin{table}[p]
\begin{center}
\caption{Medians and robust standard deviations (in parentheses) of PE, $L_2$ loss, $L_1$ loss, deviance, \#S, and FN over $100$ simulations for all methods in logistic regression by BIC and SIC, where $p = 25$}
\begin{tabular}{lllllll}
\hline
Method           &   Measures &      Lasso  &       SCAD &  MCP  &    Oracle  \\
\hline
BIC &        PE &    0.110(0.008) &    0.097(0.006) & 0.097(0.006) &    0.093(0.004)  \\
           &         $L_2$ loss &   3.055(0.656)  &     0.943(0.550) & 0.943(0.550) &    0.880(0.339) \\
           &         $L_1$ loss &   7.247(1.095)        &   1.867(1.461) & 1.867(1.461)       &   1.732(0.767)        \\
           &         Deviance &  129.36(19.20)          &   111.82(15.80) & 111.82(15.80)        &     113.12(16.05)      \\
           &         \#S &         9(2.97)  &         5(0.74) &  5(0.74) &        5(0) \\
           &         FN &         0(0)  &         0(0) &  0(0) &       0(0) \\
\hline
SIC &        PE  &   0.114(0.010)    &    0.095(0.005) & 0.095(0.005) &    0.093(0.004)  \\
           &         $L_2$ loss &     3.342(0.600)  &    0.943(0.476) & 0.943(0.476) &    0.880(0.339) \\
           &         $L_1$ loss &    7.646(1.114)       &  1.799(1.006)  & 1.799(1.006)       &  1.732(0.767)         \\
           &         Deviance &    134.93(18.35)        &   112.22(16.30)   & 112.22(16.30)     &    113.12(16.05)       \\
           &         \#S &         9(2.22)  &         5(0) &   5(0) &       5(0)  \\
           &         FN &         0(0)  &         0(0) &    0(0)  &   0(0) \\
\hline
\end{tabular}
\label{Tab1}
\end{center}
\end{table}

\begin{figure} \centering
\begin{center}%
\begin{tabular}
[l]{l}%
{\hspace{-0.45in}\includegraphics[scale=0.72]%
{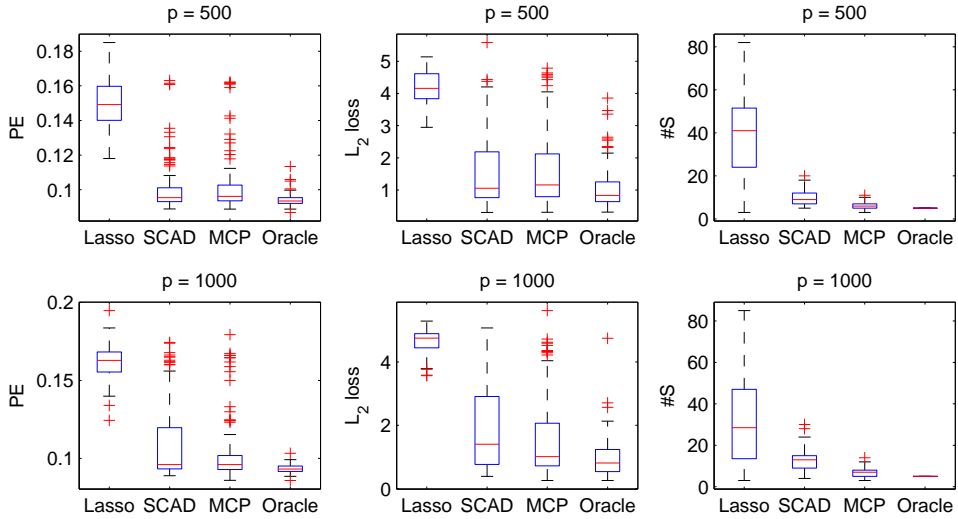}%
}%
\end{tabular}%
\vspace{-0.25in}
\caption{Boxplots of PE, $L_2$ loss, and \#S over $100$ simulations for all methods in logistic regression, where $p = 500$ and $1000$. The $x$-axis represents different methods. Top panel is for $p = 500$ and bottom panel is for $p = 1000$.}
\label{Fig2}%
\end{center}%
\end{figure}%

\begin{table}[p]
\begin{center}
\caption{Medians and robust standard deviations (in parentheses) of PE, $L_2$ loss, $L_1$ loss, deviance, \#S, and FN over $100$ simulations for all methods in logistic regression, where $p = 500$ and $1000$}
\begin{tabular}{lllllll}
\hline
$p$ & Measures &      Lasso  &       SCAD &   MCP  &    Oracle  \\
\hline
500 & PE &    0.0149(0.015) &    0.095(0.006) & 0.096(0.007) &    0.094(0.002)  \\
    & $L_2$ loss &    4.158(0.574)  &     1.054(1.054) &  1.160(0.985) &     0.834(0.452) \\
    & $L_1$ loss &      11.540(0.841)     &   2.508(2.044) &  2.481(2.292)      &   1.591(0.939)        \\
    & Deviance &        113.84(43.76)    &    100.22(16.03)   &  102.96(15.36)    &   108.06(17.33)        \\
    & \#S &        41(20.39)  &        9(3.71) & 6(1.48) &       5(0) \\
    & FN &         0(0.74)  &         0(0) & 0(0)  &        0(0) \\
\hline
1000 & PE &    0.163(0.010) &    0.096(0.020) & 0.096(0.007)  & 0.093(0.003)  \\
     & $L_2$ loss &      4.753(0.333)  &   1.400(1.591) & 1.010(1.000) &  0.808(0.517) \\
     & $L_1$ loss &      11.759(0.801)     &  3.133(3.297)  &   2.322(2.145)     &  1.490(0.949)         \\
     & Deviance &    152.19(49.36)        &   99.18(19.80)  & 103.25(16.99)      &   110.03(14.49)        \\
     & \#S &        28.5(24.83)  &     13(4.45) &  7(2.22) &       5(0) \\
     & FN &         1(0.74)  &         0(0) &    0(0)  &    0(0) \\
\hline
\end{tabular}
\label{Tab2}
\end{center}
\end{table}

\subsection{Logistic regression} \label{Sec6.1}
In this example, we demonstrate the performance of non-concave
penalized likelihood methods in logistic regression. The data were generated from the logistic regression model
(\ref{001}). We set $(n, p) = (200, 25)$ and chose the true
regression coefficients vector $\bbeta_0$ by setting
$\bbeta_1 = (2.5, -1.9, 2.8, -2.2, 3)\t$. The number of simulations was 100.
For each simulated data set, the rows of $\bX$ were sampled as
i.i.d. copies from $N(\bzero, \Sigma_0)$ with $\Sigma_0 = (0.5^{|i -
j|})_{i, j = 1, \cdots, p}$, and the response vector $\by$ was
generated independently from the Bernoulli distribution with conditional success probability vector $g(\bX \bbeta_0)$, where $g(x) =
e^x/(1+e^x)$. We compared Lasso ($L_1$ penalty), SCAD and MCP
with the oracle estimator, all of which were implemented by the ICA
algorithm to produce the solution paths. The regularization
parameter $\lambda$ was selected by BIC and the semi-Bayesian
information criterion (SIC) introduced by Lv and Liu (2008).

Six performance measures were used to compare the methods. The first measure is the prediction error (PE) defined as $E [Y - g(\bX\t \hbbeta)]^2$, where $\hbbeta$ is the estimated coefficients vector by a method and $(\bX\t, Y)$ is an independent test point. The second and third measures are the $L_2$ loss $\|\hbbeta - \bbeta_0\|_2$ and $L_1$ loss $\|\hbbeta - \bbeta_0\|_1$. The fourth measure is the deviance of the fitted model. The fifth measure, \#S, is the number of selected variables in the final model by a method in a simulation. The sixth one, FN, measures the number of missed true variables by a method in a simulation.

In the calculation of PE, an independent test sample of size 10,000 was generated to compute the expectation. For both BIC and SIC, Lasso had median FN $= 0$ with some nonzeros, and SCAD and MCP had FN $= 0$ over 100 simulations. Table \ref{Tab1} and Figure \ref{Fig1} summarize the comparison results given by PE, $L_2$ loss, $L_1$ loss, deviance, \#S, and FN, respectively for BIC and SIC.   The Lasso selects larger model sizes than SCAD and MCP.  Its associated median losses are also larger.

We also examined the performance of non-concave penalized likelihood methods in high dimensional logistic regression. The setting of this simulation is the same as above, except that $p = 500$ and $1000$. Since $p$ is larger than $n$, the information criteria break down in the tuning of $\lambda$ due to the overfitting. Thus we used five-fold cross-validation based on prediction error to select the tuning parameter. Lasso had many nonzeros of FN, and SCAD and MCP had FN $= 0$ over almost all 100 simulations except very few nonzeros. Table \ref{Tab2} and Figure \ref{Fig2} report the comparison results given by PE, $L_2$ loss, $L_1$ loss, deviance, \#S, and FN.

It is clear from Table \ref{Tab2} that LASSO selects far larger model size than SCAD and MCP.  This is due to the bias of the $L_1$ penalty.  The larger bias in LASSO forces the cross-validation to choose a smaller value of $\lambda$ to reduce its contribution to PE.  But, a smaller value of $\lambda$ allows more false positive variables to be selected.  The problem is certainly less severe for the SCAD penalty and MCP.  The performance between SCAD and MCP is comparable, as expected.

\begin{table}[htp]
\begin{center}
\caption{Medians and robust standard deviations (in parentheses) of PE, $L_2$ loss, $L_1$ loss, deviance, \#S, and FN over $100$ simulations for all methods in Poisson regression, where $p = 25$}
\begin{tabular}{llllll}
\hline
Measures &      Lasso  &       SCAD &   MCP &    Oracle  \\
\hline
PE &    7.195(2.428) &    4.081(0.826) & 4.012(0.791) &  3.688(0.574)  \\
$L_2$ loss &        0.269(0.076)  &       0.141(0.045) & 0.136(0.040) &     0.111(0.035) \\
$L_1$ loss &     0.606(0.215)      &    0.276(0.103)    & 0.271(0.094)   &   0.216(0.067)        \\
Deviance &     191.09(14.62)       &   186.73(12.72)     & 187.23(13.14)   &    187.72(15.28)       \\
\#S &         9(2.22)  &         5(0.74) &    5(0.74)  &   5(0) \\
FN &         0(0)  &         0(0) &    0(0)   &  0(0) \\
\hline
\end{tabular}
\label{Tab3}
\end{center}
\end{table}

\begin{table}[hp]
\begin{center}
\caption{Medians and robust standard deviations (in parentheses) of PE, $L_2$ loss, $L_1$ loss, deviance, \#S, and FN over $100$ simulations for all methods in Poisson regression by BIC and CV, where $p = 500$ and $1000$}
\begin{tabular}{llllllll}
\hline
$p$ & Method           &   Measures &      Lasso  &       SCAD &  MCP  &    Oracle  \\
\hline
500 & BIC &    PE &    26.989(11.339) &    4.820(1.772) & 4.672(1.593)  &  3.479(0.738)  \\
 &      &     $L_2$ loss &  0.790(0.206)  &  0.199(0.074) & 0.178(0.076) &  0.104(0.043) \\
 &   &    $L_1$ loss &  2.446(0.638)         &  0.424(0.165)    &   0.371(0.161)  &     0.184(0.083)      \\
     &      &         Deviance & 202.65(22.23)          &     187.15(16.41)  &  189.66(17.24)  &    189.30(21.73)       \\
     &      &         \#S &         29(5.93)  &         9(3.34) &  7(2.22)     &  5(0) \\
     &      &         FN &         0(0)  &         0(0) &   0(0)  &    0(0) \\
\cline{2-7}
     & CV &        PE  &   24.200(9.636)    &   4.542(1.554) & 4.272(1.503)   & 3.479(0.738)  \\
     &      &         $L_2$ loss &      0.698(0.162)  &  0.168(0.065) & 0.168(0.057) & 0.104(0.043) \\
     &      &         $L_1$ loss &   3.229(1.368)        &   0.495(0.201)        & 0.411(0.162) &  0.184(0.083)        \\
     &      &         Deviance &  117.16(40.13)          &  166.58(21.76)    & 173.01(19.17) &  189.30(21.73)        \\
     &      &         \#S &         63.5(24.83)  &         18(10.75) &   12.5(6.67)  &    5(0)  \\
     &      &         FN &         0(0)  &         0(0) &    0(0)  &   0(0) \\
\hline
1000 & BIC &        PE &    33.069(14.089) &    5.523(2.027) & 5.144(1.808) &  3.676(0.772)  \\
     &      &         $L_2$ loss &    0.971(0.209)  &    0.210(0.094) &  0.187(0.088)  &  0.108(0.047) \\
     &      &         $L_1$ loss &   2.990(0.689)        &   0.485(0.232)    & 0.443(0.198)   &  0.197(0.090)         \\
     &      &         Deviance &     199.99(22.89)       &    180.34(13.07)  & 181.21(15.31)   &   187.98(17.22)       \\
     &      &         \#S &         34(7.41)  &        11.5(4.08) & 9(2.22)    &    5(0) \\
     &      &         FN &         0(0)  &         0(0) &   0(0)   &   0(0) \\
\cline{2-7}
     & CV &        PE  &   31.701(16.571)    &  4.821(1.732) & 4.700(1.702) &  3.676(0.772)  \\
     &      &         $L_2$ loss &    0.889(0.201)  &   0.162(0.077) &  0.162(0.064) &  0.108(0.047) \\
     &      &         $L_1$ loss &   4.297(1.646)      &  0.506(0.341)  &   0.454(0.239)  &    0.197(0.090)       \\
     &      &         Deviance &   92.89(44.51)         &   160.23(20.80) &   169.34(23.44)  &   187.98(17.22)        \\
     &      &         \#S &  83(40.77)  &    22(11.86) &  14(7.04)    &   5(0)  \\
     &      &         FN &         0(0)  &         0(0) &   0(0)    &  0(0) \\
\hline
\end{tabular}
\label{Tab4}
\end{center}
\end{table}

\begin{table}[htp]
\begin{center}
\caption{Classification errors in the neuroblastoma data set}
\begin{tabular}{lccccc}
\hline
& \multicolumn{2}{c}{3-year EFS}  &   &  \multicolumn{2}{c}{Gender}  \\
\cline{2-3}  \cline{5-6}
Method          &    \# of genes  &   Test error    &      &    \# of genes  &   Test error  \\
\hline
Lasso &        56 &    23/114 &   &  4 & 5/126 \\
SCAD &        10  &   18/114  &    &  2 & 4/126 \\
MCP &       7   &   23/114  &    &  1 & 12/126 \\
SIS &        5  &   19/114  &    &  6 & 4/126 \\
ISIS &        23  &   22/114  &    &  2 & 4/126 \\
\hline
\end{tabular}
\label{Tab5}
\end{center}
\end{table}

\begin{table}[hp]
\begin{center}
\caption{Selected genes for the 3-year EFS classification}
\begin{tabular}{lccclccc}
\hline
Gene  &  Lasso &  SCAD & MCP & Gene  &  Lasso &  SCAD & MCP  \\
\hline
\verb'A_24_P182182' &  & x &   & \verb'Hs419768.1' & x &  & \\
\verb'A_23_P144096' & x & & & \verb'A_23_P313728' & x &   & \\
\verb'A_23_P124514' & x & & & \verb'A_23_P12884' & x & & \\
\verb'A_23_P502879' & x & & & \verb'A_23_P130626' & x & & \\
\verb'A_23_P71319' & x & & & \verb'A_23_P143958' &  & x & \\
\verb'A_24_P73158' & x & x & & \verb'Hs155462.1' &  & x & \\
\verb'A_24_P282251' & x & & x & \verb'A_23_P209394' & x & &  \\
\verb'A_23_P125435' & x &  & & \verb'A_24_P100419' & x &  & \\
\verb'A_23_P80491' & x &  & & \verb'Hs379382.1' & x &  & \\
\verb'A_23_P77779' & x &  & & \verb'A_24_P271696' & x & & \\
\verb'A_23_P36076' & x & x & & \verb'Hs381187.1' & x &  & \\
\verb'A_23_P35349' & x &  & & \verb'Hs265827.1' & x &  & \\
\verb'A_23_P208030' & x &  & & \verb'Hs269914.3' & x &  & \\
\verb'A_23_P72737' & x &  & x & \verb'Hs36034.1' & x & & \\
\verb'A_23_P25194' & x &  & & \verb'A_23_P83751' & x &  & \\
\verb'A_23_P200043' & x &  & & \verb'A_23_P96325' & x & & \\
\verb'A_23_P422809' & x &  & & \verb'A_23_P97553' & x & & \\
\verb'A_23_P110345' & x &  & x & \verb'A_24_P232158' & x & & \\
\verb'A_23_P5131' & x &  & x & \verb'A_23_P9836' & x & & \\
\verb'A_23_P11859' & x &  & & \verb'Hs170298.1' & x &  & x \\
\verb'A_23_P7376' & x &  & & \verb'r60_a135' & x &  & \\
\verb'A_23_P211985' & x &  & & \verb'A_23_P133956' & x & & \\
\verb'A_24_P365954' &  & x & & \verb'A_32_P27511' & x &  & \\
\verb'A_23_P86975' & x &  & & \verb'A_23_P80626' & x & & \\
\verb'A_23_P89910' & x &  & & \verb'A_32_P158708' & x & & \\
\verb'A_24_P285055' & x &  & & \verb'A_23_P100764' &  & x & \\
\verb'A_23_P68547' & x &  & & \verb'Hs407755.1' & x &  & \\
\verb'A_23_P6252' & x &  & & \verb'Hs86643.1' & x & & \\
\verb'A_23_P386356' & x &  & & \verb'Hs422789.1' & x & x & \\
\verb'A_24_P50228' & & & x & \verb'A_23_P131866' &  & x & \\
\verb'Hs37637.1' & x &  & & \verb'A_23_P147397' & x &  & \\
\verb'Hs457415.1' &  & x & x & \verb'A_23_P13852' & x & & \\
\hline
\end{tabular}
\label{Tab6}
\end{center}
\end{table}

\begin{table}[htp]
\begin{center}
\caption{Selected genes for the gender classification}
\begin{tabular}{lccc}
\hline
Gene  &  Lasso &  SCAD & MCP \\
\hline
\verb'A_23_P329835' &  x   &  & \\
\verb'A_23_P259314' &  x   &   &    \\
\verb'A_23_P137238' &  x   & x & x \\
\verb'A_24_P500584' & x   & x & \\
\hline
\end{tabular}
\label{Tab7}
\end{center}
\end{table}

\subsection{Poisson regression} \label{Sec6.2}
In this example, we demonstrate the performance of non-concave penalized likelihood methods in Poisson regression. The data were generated from the Poisson regression model (\ref{001}). The setting of this example is similar to that in Section \ref{Sec6.1}. We set $(n, p) = (200, 25)$ and chose the true regression coefficients vector $\bbeta_0$ by setting $\bbeta_1 = (1.25, -0.95, 0.9, -1.1, 0.6)\t$. For each simulated data set, the response vector $\by$ was generated independently from the Poisson distribution with conditional mean vector $\exp(\bX \bbeta_0)$. The regularization parameter $\lambda$ was selected by BIC (SIC performed similarly to BIC).

The PE is defined as $E [Y - \exp(\bX\t \hbbeta)]^2$, where $\hbbeta$ is the estimated coefficients vector by a method and $(\bX\t, Y)$ is an independent test point. Lasso, SCAD and MCP had FN $= 0$ over 100 simulations. Table \ref{Tab3} summarizes the comparison results given by PE, $L_2$ loss, $L_1$ loss, deviance, \#S, and FN.

We also examined the performance of non-concave penalized likelihood methods in high dimensional Poisson regression. The setting of this simulation is the same as above, except that $p = 500$ and $1000$. The regularization parameter $\lambda$ was selected by BIC and five-fold cross-validation (CV) based on prediction error. For both BIC and CV, Lasso had median FN $= 0$ with some nonzeros, and SCAD and MCP had FN $= 0$ over 100 simulations. Table \ref{Tab4} reports the comparison results given by PE, $L_2$ loss, $L_1$ loss, deviance, \#S, and FN.

\subsection{Real data analysis} \label{Sec6.3}
In this example, we apply non-concave penalized likelihood methods to the neuroblastoma data set, which was studied by Oberthuer \etal (2006). This data set, obtained via the MicroArray Quality Control phase-II (MAQC-II) project, consists of gene expression profiles for 10,707 genes from 251 patients of the German Neuroblastoma Trials NB90-NB2004, diagnosed between 1989 and 2004. The patients at diagnosis were aged from 0 to 296 months with a median age of 15 months. The study aimed to develop a gene expression-based classifier for neuroblastoma patients that can reliably predict courses of the disease.

We analyzed this data set for two binary responses: 3-year event-free survival (3-year EFS) and gender, where 3-year EFS indicates whether a patient survived 3 years after the diagnosis of neuroblastoma.  There are 246 subjects with 101 females and 145 males, and 239 of them have the 3-year EFS information available (49 positives and 190 negatives). We applied Lasso, SCAD and MCP using the logistic regression model. Five-fold cross-validation was used to select the tuning parameter. For the 3-year EFS classification, we randomly selected 125 subjects (25 positives and 100 negatives) as the training set and the rest as the test set. For the gender classification, we randomly chose 120 subjects (50 females and 70 males) as the training set and the rest as the test set. Table \ref{Tab5} reports the classification results of all methods, as well as those of SIS and ISIS, which were extracted from Fan \etal (2009). Tables \ref{Tab6} and \ref{Tab7} list the selected genes by Lasso, SCAD and MCP for the 3-year EFS classification and gender classification, respectively.

\section{Discussions} \label{Sec7}
We have studied penalized likelihood methods for ultra-high
dimensional variable selection. In the context of GLMs, we have shown that such methods have model selection consistency with oracle properties even
for NP-dimensionality, for a class of non-concave penalized likelihood approaches. Our results are consistent with a known fact in the literature that concave penalties can reduce the bias problems of convex penalties. The convex function of $L_1$-penalty falls at the boundary of the class of penalty functions under consideration. We have used the coordinate
optimization to find the solution paths and illustrated the performance of non-concave penalized likelihood methods with numerical studies. Our results show that the coordinate
optimization works equally well and efficiently for producing the entire solution paths for concave penalties.

\section{Proofs} \label{Sec8}

\subsection{Proof of Theorem \ref{T1}} \label{Sec8.1}
We will first derive the necessary condition.  In view of (\ref{002}), we have
\begin{equation} \label{123}
\nabla \ell_n(\bbeta) =  n^{-1} \left[\bX\t \by - \bX\t
\bmu\left(\btheta\right)\right] \ \text{ and } \  \nabla^2
\ell_n(\bbeta) = - n^{-1} \bX\t \Sig\left(\btheta\right) \bX,
\end{equation}
where $\btheta = \bX \bbeta$. It follows from the classical
optimization theory that if $\hbbeta = (\hbeta_1, \cdots, \hbeta_p)\t$ is a
local maximizer of the penalized likelihood (\ref{004}), it satisfies the Karush-Kuhn-Tucker (KKT) conditions, i.e., there exists some $\bv =
(v_1, \cdots, v_p)\t \in \mathbf{R}^p$ such that
\begin{equation} \label{078}
\bX\t \by - \bX\t \bmu(\hbtheta) -  n \lambda_n \bv = \bzero,
\end{equation}
where $\hbtheta = \bX \hbbeta$, $v_j = \bar{\rho}(\hbeta_j)$ for
$\hbeta_j \neq 0$, and $v_j \in [-\rho'(0+), \rho'(0+)]$ for
$\hbeta_j = 0$. Let $\calS = \supp(\hbbeta)$.  Note that $\hbbeta$ is also a local maximizer of
(\ref{004}) constrained on the $\|\hbbeta\|_0$-dimensional subspace
$\mathcal{B} = \{\bbeta \in \mathbf{R}^p: \bbeta_c = \bzero\}$ of
$\mathbf{R}^p$, where $\bbeta_c$ denotes the subvector of $\bbeta$
formed by components in $\calS^c$, the complement of
$\calS$. It follows from the second order condition that
\begin{equation} \label{102}
\lambda_{\min} \left[\bX_1\t \Sig\left(\hbtheta\right) \bX_1\right]
\geq  n \lambda_n \kappa(\rho; \hbbeta_1),
\end{equation}
where  $\kappa(\rho; \hbbeta_1)$ is given by (\ref{016}).
It is easy to see that equation (\ref{078}) can be equivalently written as
\begin{align}
\label{079}
& \bX_1\t \by - \bX_1\t \bmu(\hbtheta) -  n \lambda_n \bar{\rho}(\hbbeta_1) = \bzero, \\
\label{080}
& \|\bz\|_\infty \leq \rho'(0+),
\end{align}
where $\bz = ( n \lambda_n)^{-1} \bX_2\t [\by - \bmu(\hbtheta)]$ and $\bX_2$ denotes the submatrix of $\bX$ formed by columns in $\calS^c$.

We now prove the sufficient condition. We first constrain the penalized likelihood (\ref{004}) on the
$\|\hbbeta\|_0$-dimensional subspace $\mathcal{B}$ of
$\mathbf{R}^p$. It follows from
condition (\ref{023}) that  $Q_n(\bbeta)$ is strictly concave in a ball $\mathcal{N}_0$ in the
subspace $\mathcal{B}$ centered at $\hbbeta$. This along with
equation (\ref{021}) immediately entails that $\hbbeta$, as a
critical point of $Q_n(\bbeta)$ in $\mathcal{B}$, is the unique
maximizer of $Q_n(\bbeta)$ in the neighborhood $\mathcal{N}_0$.

It remains to prove that the sparse vector $\hbbeta$ is indeed a
strict local maximizer of $Q_n(\bbeta)$ on the space
$\mathbf{R}^p$. To show this, take a sufficiently small ball $\mathcal{N}_1$ in $\mathbf{R}^p$ centered at $\hbbeta$  such that $\mathcal{N}_1 \cap \mathcal{B}
\subset \mathcal{N}_0$. We then need to show that
$Q_n(\hbbeta) > Q_n(\bgamma_1)$ for any $\bgamma_1 \in
\mathcal{N}_1 \setminus \mathcal{N}_0$.  Let $\bgamma_2$ be the projection
of $\bgamma_1$ onto the subspace $\mathcal{B}$.  Then we have $\bgamma_2
\in \mathcal{N}_0$, which entails that
$Q_n(\hbbeta) > Q_n(\bgamma_2)$ if $\bgamma_2 \neq \hbbeta$, since $\hbbeta$ is the strict maximizer of $Q_n(\bbeta)$ in $\mathcal{N}_0$. Thus, it suffices to show that $Q_n(\bgamma_2) > Q_n(\bgamma_1)$.

By the mean-value theorem, we have
\begin{equation} \label{082}
Q_n(\bgamma_1) - Q_n(\bgamma_2) = \nabla\t Q_n(\bgamma_0) (\bgamma_1 - \bgamma_2),
\end{equation}
where $\bgamma_0$ lies on the line segment joining $\bgamma_2$ and $\bgamma_1$.
Note that the components of $\bgamma_1 - \bgamma_2$ are zero for the indices in $\calS$ and the sign of $\gamma_{0,j}$ is the same as that of $\gamma_{1, j}$ for $j \not \in \calS$, where $\gamma_{0,j}$ and $\gamma_{1,j}$ are the $j$-th components of $\bgamma_0$ and $\bgamma_1$, respectively. Therefore, the right hand side of (\ref{082}) can be expressed as
\begin{eqnarray}   \label{08a}
 & & \left\{n^{-1} \bX_2\t \left[\by - \bmu(\bX \bgamma_0)\right]\right\}\t \bgamma_{1, 2}  - \lambda_n \sum_{j \not \in \calS} \rho'(|\gamma_{0, j}|) |\gamma_{1, j}|,
\end{eqnarray}
where $\bgamma_{1, 2}$ is a subvector of $\bgamma_1$ formed by the
components in $\calS^c$. By $\bgamma_1 \in
\mathcal{N}_1 \setminus \mathcal{N}_0$, we have $\bgamma_{1, 2} \neq \bzero$.

It follows from the concavity of $\rho$ in Condition \ref{con1} that
$\rho'(t)$ is decreasing in $t \in [0, \infty)$. By condition (\ref{022}) and the continuity of $\rho'(t)$ and $b'(\theta)$,
there exists some $\delta > 0$ such that
for any $\bbeta$ in a ball in $\bR^p$ centered at $\hbbeta$ with radius $\delta$,
\begin{equation} \label{081}
\| ( n \lambda_n)^{-1}  \bX_2\t \left[\by - \bmu(\bX \bbeta)\right]
 \|_\infty < \rho'(\delta).
\end{equation}
We further shrink the radius of the ball $\mathcal{N}_1$ to less
than $\delta$ so that $|\gamma_{0, j}| \leq |\gamma_{1,j}| < \delta$
for $j \not \in \calS$ and (\ref{081}) holds for any $\bbeta \in
\mathcal{N}_1$.  Since $\bgamma_0 \in \mathcal{N}_1$, it follows
from (\ref{081}) that the term (\ref{08a}) is strictly
less than
$$
\lambda_n  \rho'(\delta) \|\bgamma_{1, 2}\|_1 -  \lambda_n  \rho'(\delta) \|\bgamma_{1, 2}\|_1 = 0,
$$
where the monotonicity of $\rho'(\cdot)$ was used in the second term.  Thus we conclude that $Q_n(\gamma_1) < Q_n(\gamma_2)$.  This completes the proof.

\subsection{Proof of Proposition \ref{P2}} \label{Sec8.2}
Let $ \partial \mathcal{L}_c = \left\{\bbeta \in \mathbf{R}^p: \ell_n(\bbeta) = c\right\}$ be the level set.
By the concavity of $\ell_n(\bbeta)$, we can easily show that for $c < \ell_n(\bzero)$, $\mathcal{L}_c$ is a closed convex set with $\bbeta_*$ and $\bzero$ being its interior points and the level set $\partial \mathcal{L}_c$ is its boundary. We now show that the global maximizer of the penalized likelihood $Q_n(\bbeta)$ belongs to $\mathcal{L}_c$.

For any $\bgamma \in \partial \mathcal{L}_c$, let $\Gamma_{\bgamma} = \{t \bgamma: t \in (1, \infty)\}$ be a ray. By the convexity of $\mathcal{L}_c$, we have $\{t \bgamma: t \in [0, 1]\} \subset \mathcal{L}_c$ for $\bgamma \in \partial \mathcal{L}_c$, which implies that
\[ \bigcup_{\bgamma \in \partial \mathcal{L}_c} \Gamma_{\bgamma} = \mathbf{R}^p \setminus \mathcal{L}_c. \]
Thus to show that the global maximizer of $Q_n(\bbeta)$ belongs to
$\mathcal{L}_c$, it suffices to prove $Q_n(t \bgamma) <
Q_n(\bgamma)$ for any $t \in (1, \infty)$ and $\bgamma \in \partial
\mathcal{L}_c$. This follows easily from the definition of $Q_n(\bbeta)$,
$\ell_n(t \bgamma) < c = \ell_n(\bgamma)$, and $\sum_{j = 1}^p
p_{\lambda_n}(t |\gamma_j|) \geq \sum_{j = 1}^p
p_{\lambda_n}(|\gamma_j|)$, where $\bgamma = (\gamma_1, \cdots,
\gamma_p)\t$.

It remains to prove that the local maximizer of $Q_n(\bbeta)$ in $\mathcal{L}_c$ must be a global maximizer. This is entailed by the concavity of $Q_n(\bbeta)$ on $\mathcal{L}_c$, which is ensured by condition (\ref{141}). This concludes the proof.

\subsection{Proof of Proposition \ref{P3}} \label{Sec8.3}
Since $c < \ell_n(\bzero)$, from the proof of Proposition \ref{P2}
we know that the global maximizer of the penalized likelihood
$Q_n(\bbeta)$ belongs to $\mathcal{L}_c$. Note that by assumption, the SCAD penalized likelihood estimator $\hbbeta = (\hbeta_1, \cdots, \hbeta_p)\t \in \mathcal{L}_c$ and $\min_{j = 1}^p |\hbeta_j| > a \lambda_n$. It follows from (\ref{004}) and (\ref{003}) that $\hbbeta$ is a critical point of $\ell_n(\bbeta)$ and thus $\hbbeta = \bbeta_*$ by the strict concavity of $\ell_n(\bbeta)$. It remains to prove that $\bbeta_*$ is the maximizer of $Q_n(\bbeta)$ on $\mathcal{L}_c$.

The key idea is to use a first order Taylor expansion of $\ell_n(\bbeta)$ around $\bbeta_*$
and retain the Lagrange remainder term. This along with $\nabla
\ell_n(\bbeta_*) = \bzero$ and $\min_{\bbeta \in \mathcal{L}_c}
\lambda_{\min} [n^{-1} \bX\t \Sig(\bX \bbeta) \bX] \geq c_0 $ gives
for any $\bbeta \in \mathcal{L}_c$,
\[
Q_n(\bbeta) \leq \widetilde{Q}_n(\bbeta) \equiv \ell_n(\bbeta_*) -\frac{c_0}{2} \left\|\bbeta - \bbeta_*\right\|_2^2 -  \sum_{j = 1}^p p_{\lambda_n}(|\beta_j|),
\]
since $\bbeta_*$ is in the convex set $\mathcal{L}_c$.
Thus if $\bbeta_*$ is the global maximizer of $\widetilde{Q}_n(\bbeta)$ on $\mathbf{R}^p$, then we have for any $\bbeta \in \mathcal{L}_c$,
$$
Q_n(\bbeta) \leq \widetilde{Q}_n(\bbeta) \leq \widetilde{Q}_n(\bbeta_*) =
{Q}_n(\bbeta_*).
$$
This entails that $\bbeta_*$ is the global maximizer of $Q_n(\bbeta)$.

To maximize $\widetilde{Q}_n(\bbeta)$, we only need to maximize it
componentwise. Let $\bbeta_* = (\beta_{*, 1}, \cdots, \beta_{*,
p})\t$. Then it remains to show that for each $j = 1, \cdots, p$,
$\beta_{*, j}$ is the global minimizer of the univariate SCAD
penalized least squares problem
\begin{equation} \label{142}
\min\limits_{\beta \in \mathbf{R}} g_j(\beta) = \min\limits_{\beta \in \mathbf{R}} \left\{\frac{c_0}{2} \left(\beta - \beta_{*, j}\right)^2 + p_{\lambda_n}(|\beta|)\right\}.
\end{equation}
This can easily been shown from the analytical solution to
(\ref{142}). For the sake of completeness, we give a simple proof here.

Recall that we have shown that $\hbbeta = \bbeta_*$. In view of (\ref{142}) and $|\beta_{*, j}| > a \lambda_n$, for any
$|\beta| > a \lambda_n$ with $\beta \neq \beta_{*, j}$, we have
\[
  g_j(\beta) > p_{\lambda_n}(|\beta|) =  p_{\lambda_n}(a \lambda_n) = g_j(|\beta_{*, j}|),
\]
where we used the fact that $p_{\lambda_n}(\cdot)$ is constant on $[a \lambda_n, \infty)$.  Thus, it suffices to prove $g_j(\beta) > p_{\lambda_n}(a \lambda_n)$ on the interval $|\beta|\leq a \lambda_n$.   For such a $\beta$, we have $p_{\lambda_n}(a \lambda_n) - p_{\lambda_n}(|\beta|) \leq \lambda_n (a \lambda_n - |\beta|)$. Thus we need to show that
\[
\min\limits_{z \in [0, a \lambda_n]} \left\{\frac{c_0}{2} \left(|\beta_{*, j}| - a \lambda_n + z\right)^2 - \lambda_n z\right\} > 0,
\]
which always holds as long as $|\beta_{*, j}| > (a + \frac{1}{2 c_0}) \lambda_n$ and thus completes the proof.

\subsection{Proof of Proposition \ref{P4}} \label{Sec8.4}
Let $\mathcal{A}$ be any $s$-dimensional coordinate subspace different from $\mathcal{A}_1 = \{(\beta_1, \cdots, \beta_p)\t \in \mathbf{R}^p: \beta_j = 0 \text{ for } j \notin \supp(\hbbeta)\}$. Clearly $\mathcal{A}_1 \oplus \mathcal{A}$ is a $d$-dimensional coordinate subspace with $d \leq 2 s$. Then part a) follows easily from the assumptions and Proposition \ref{P2}. Part b) is an easy consequence of Proposition \ref{P3} in view of the assumptions and the fact that
\[ \max_{t \in [0, \infty)} p_{\lambda_n}(t) = p_{\lambda_n}(a \lambda_n) = \frac{(a + 1) \lambda_n^2}{2} \]
for the SCAD penalty $p_\lambda$ given by (\ref{003}).

\subsection{Proof of Proposition \ref{P1}} \label{Sec8.5}
Part a) follows easily from a simple application of Hoeffding's inequality (Hoeffding, 1963), since $a_1 Y_1, \cdots, a_n Y_n$ are $n$ independent bounded random variables, where $\ba = (a_1, \cdots, a_n)\t$. We now prove part b). In view of condition (\ref{010}), $a_i Y_i - a_i b'(\theta_{0, i})$ are $n$ independent random variables with mean zero and satisfy
\begin{align*}
E \left|a_i Y_i - a_i b'\left(\theta_{0, i}\right)\right|^m & = |a_i|^m E \left|Y_i - b'\left(\theta_{0, i}\right)\right|^m \leq |a_i|^m m! M^{m - 2} \frac{v_0}{2}\\
& \leq \frac{m!}{2} \left(\left\|\ba\right\|_\infty M\right)^{m - 2} a_i^2 v_0, \quad m \geq 2.
\end{align*}
Thus an application of Bernstein's inequality (see, e.g., Bennett, 1962 or van der Vaart and Wellner, 1996) yields
\begin{align*}
P\left(\left|\ba\t \bY - \ba\t \bmu(\btheta_0)\right| > \varepsilon\right) & \leq 2 \exp\left[-\frac{1}{2} \frac{\varepsilon^2}{\sum_{i = 1}^n a_i^2 v_0 + \left\|\ba\right\|_\infty M \varepsilon}\right] \\
& = 2 \exp\left[-\frac{1}{2} \frac{\varepsilon^2}{\left\|\ba\right\|_2^2 v_0 + \|\ba\|_\infty M \varepsilon}\right],
\end{align*}
which concludes the proof.

\subsection{Proof of Theorem \ref{T2}} \label{Sec8.6}
We break the whole proof into several steps. Let $\bX_1$ and $\bX_2$
respectively be the submatrices of $\bX$ formed by columns in
$\mathfrak{M}_0 = \supp(\bbeta_0)$ and its complement
$\mathfrak{M}_0^c$, and $\btheta_0 = \bX \bbeta_0$.  Let $\bxi =
(\xi_1, \cdots, \xi_p)\t = \bX\t \by - \bX\t \bmu(\btheta_0)$.
Consider events
\[
\mathcal{E}_1 = \left\{\left\|\bxi_{\mathfrak{M}_0}\right\|_\infty \leq c_1^{-1/2} \sqrt{n \log n}\right\} \quad \text{and} \quad \mathcal{E}_2 = \left\{\left\|\bxi_{\mathfrak{M}_0^c}\right\|_\infty \leq u_n \sqrt{n}\right\},
\]
where $u_n = c_1^{-1/2} n^{1/2 - \alpha} (\log n)^{1/2}$ is a diverging
sequence  and $\bv_A$ denotes a subvector of $\bv$ consisting of
elements in $A$. Since $\|\bx_j\|_2 = \sqrt{n}$, it follows from
Bonferroni's inequality and (\ref{137}) that
\begin{align}\label{106}
& \quad P\left(\mathcal{E}_1 \cap \mathcal{E}_2\right)  \\
& \geq 1 - \sum_{j \in \mathfrak{M}_0} P\left(|\xi_j| > c_1^{-1/2}
\sqrt{n \log n}\right) - \sum_{j \in \mathfrak{M}_0^c}
P\left(|\xi_j| > u_n \sqrt{n}\right)
\nonumber \\
& \geq 1 - 2 \left[s n^{-1} + \left(p - s\right)
e^{-c_1 u_n^2}\right] \nonumber \\
& = 1 - 2 [s n^{-1} + (p - s) e^{-n^{1 - 2 \alpha}
\log n}], \nonumber
\end{align}
where $s = \|\bbeta_0\|_0$ and $u_n \leq \sqrt{n}/\max_{j = 1}^p
\|\bx_j\|_\infty$ for unbounded responses, which is guaranteed for sufficiently large $n$ by
Condition \ref{con3}. Under the event $\mathcal{E}_1 \cap \mathcal{E}_2$, we
will show that there exists a solution $\hbbeta \in \bR^p$ to
(\ref{021})--(\ref{023}) with $\sgn(\hbbeta) = \sgn(\bbeta_0)$ and
$\|\hbbeta - \bbeta_0\|_\infty = O(n^{-\gamma} \log n)$, where the
function $\sgn$ is applied componentwise.

\smallskip

\textit{Step 1: Existence of a solution to equation (\ref{021})}. We first prove that for  sufficiently large $n$, equation (\ref{021}) has a solution $\hbbeta_1$ inside the hypercube
\[
\mathcal{N} = \left\{\bdelta \in \mathbf{R}^s: \|\bdelta - \bbeta_1\|_\infty = n^{-\gamma} \log n\right\}.
\]
For any $\bdelta = (\delta_1, \cdots, \delta_s)\t \in
\mathcal{N}$, since $d_n \geq n^{-\gamma} \log n$, we have
\begin{equation} \label{107}
\min\nolimits_{j = 1}^s |\delta_j| \geq \min\nolimits_{j \in \mathfrak{M}_0} |\beta_{0, j}| - d_n = d_n
\end{equation}
and $\sgn(\bdelta) = \sgn(\bbeta_1)$. Let $\bet =  n \lambda_n
\bar{\rho}(\bdelta)$.   Using the monotonicity condition of
$\rho'(t)$, by (\ref{107}) we have
\[
\|\bet\|_\infty \leq  n \lambda_n \rho'(d_n),
\]
which along with the definition of $\mathcal{E}_1$ entails
\begin{equation} \label{109}
\|\bxi_{\mathfrak{M}_0} - \bet\|_\infty \leq c_1^{-1/2} \sqrt{n \log
n} +  n \lambda_n \rho'(d_n).
\end{equation}

Define vector-valued functions
\[ \bgamma(\bdelta) = (\gamma_1(\bdelta), \cdots, \gamma_p(\bdelta))\t = \bX\t \bmu(\bX_1 \bdelta), \quad \bdelta \in \mathbf{R}^s \]
and
\[ \bPsi(\bdelta) = \bgamma_{\mathfrak{M}_0}(\bdelta) - \bgamma_{\mathfrak{M}_0}(\bbeta_1) - (\bxi_{\mathfrak{M}_0} - \bet), \quad \bdelta \in \mathbf{R}^s. \]
Then, equation (\ref{021}) is equivalent to $\bPsi(\bdelta) = \bzero$.  We need to show that the latter has a solution inside the hypercube $\mathcal{N}$. To this end, we represent $\bgamma_{\mathfrak{M}_0}(\bdelta)$ by using a second order Taylor expansion around $\bbeta_1$ with the Lagrange remainder term componentwise and obtain
\begin{equation} \label{110}
\bgamma_{\mathfrak{M}_0}(\bdelta) =
\bgamma_{\mathfrak{M}_0}(\bbeta_1) +  \bX_1\t
\Sig\left(\btheta_0\right) \bX_1 (\bdelta - \bbeta_1) + \br,
\end{equation}
where $\br = (r_1, \cdots, r_s)\t$ and for each $j = 1, \cdots, s$,
\[ r_j = \frac{1}{2} \left(\bdelta - \bbeta_1\right)\t \nabla^2 \gamma_{j}(\bdelta_j) \left(\bdelta - \bbeta_1\right) \]
with $\bdelta_j$ some $s$-vector lying on the line segment joining $\bdelta$ and $\bbeta_1$. By (\ref{009}), we have
\begin{align} \label{111}
\left\|\br\right\|_\infty & \leq \max_{\bdelta_0 \in \mathcal{N}} \max_{j = 1}^s \frac{1}{2} \lambda_{\max}\left[\bX_1\t \diag\left\{\left|\bx_j\right| \circ \left|\bmu''\left(\bX_1 \bdelta_0\right)\right|\right\} \bX_1\right] \left\|\bdelta - \bbeta_1\right\|_2^2 \\
\nonumber
&  = O\left[s n^{1 - 2 \gamma} (\log n)^2\right].
\end{align}
Let
\begin{equation} \label{112}
\overline{\bPsi}(\bdelta) \equiv  \left[\bX_1\t
\Sig\left(\btheta_0\right) \bX_1\right]^{-1} \bPsi(\bdelta) =
\bdelta - \bbeta_1 + \bu,
\end{equation}
where  $\bu = -[\bX_1\t \Sig(\btheta_0) \bX_1]^{-1}
(\bxi_{\mathfrak{M}_0} - \bet - \br)$. It follows from (\ref{109}),
(\ref{111}), and (\ref{007}) in Condition \ref{con2} that for any
$\bdelta \in \mathcal{N}$,
\begin{align} \label{140}
\left\|\bu\right\|_\infty
& \leq  \left\|\left[\bX_1\t \Sig\left(\btheta_0\right) \bX_1\right]^{-1}\right\|_\infty \left(\|\bxi_{\mathfrak{M}_0} - \bet\|_\infty + \|\br\|_\infty\right)\\
\nonumber
& = O\left[b_s n^{-1/2} \sqrt{\log n} + b_s \lambda_n \rho'(d_n) + b_s s n^{- 2 \gamma} (\log n)^2\right].
\end{align}
By Condition \ref{con3}, the first and third terms are of order
$o(n^{-\gamma} \log n)$ and so is the second term by
(\ref{104}). This shows that
\[
   \left\|\bu\right\|_\infty  = o(n^{-\gamma} \log n). \]
By (\ref{112}),  for sufficiently large $n$, if $(\bdelta -
\bbeta_1)_j = n^{-\gamma} \sqrt{\log n}$, we have
\begin{equation} \label{114}
\overline{\Psi}_j(\bdelta) \geq n^{-\gamma} \sqrt{\log n} - \left\|\bu\right\|_\infty  \geq 0,
\end{equation}
and if $(\bdelta - \bbeta_1)_j = -n^{-\gamma} \sqrt{\log n}$, we have
\begin{equation} \label{115}
\overline{\Psi}_j(\bdelta) \leq -n^{-\gamma} \sqrt{\log n} + \left\|\bu\right\|_\infty \leq 0,
\end{equation}
where $\overline{\bPsi}(\bdelta) = (\overline{\Psi}_1(\bdelta), \cdots, \overline{\Psi}_s(\bdelta))\t$. By the continuity of the vector-valued function $\overline{\bPsi}(\bdelta)$, (\ref{114}) and (\ref{115}), an application of Miranda's existence theorem (see, e.g., Vrahatis, 1989) shows that equation $\overline{\bPsi}(\bdelta) = \bzero$ has a solution $\hbbeta_1$ in $\mathcal{N}$. Clearly $\hbbeta_1$ also solves equation $\bPsi(\bdelta) = \bzero$ in view of (\ref{112}). Thus we have shown that equation (\ref{021}) indeed has a solution $\hbbeta_1$ in $\mathcal{N}$.

\smallskip

\textit{Step 2: Verification of condition (\ref{022})}.  Let
$\hbbeta \in \mathbf{R}^p$ with $\hbbeta_{\mathfrak{M}_0} =
\hbbeta_1 \in \mathcal{N}$ a solution to equation (\ref{021}) and
$\hbbeta_{\mathfrak{M}_0^c} = \bzero$, and $\hbtheta = \bX \hbbeta$.
We now show that $\hbbeta$ satisfies inequality (\ref{022}) for
$\lambda_n$ given by (\ref{104}). Note that
\begin{align} \label{116}
\bz & = \left( n \lambda_n\right)^{-1} \left\{\left[\bX_2\t \by - \bX_2\t \bmu\left(\btheta_0\right)\right] - \left[\bX_2\t \bmu\left(\hbtheta\right) - \bX_2\t \bmu\left(\btheta_0\right)\right]\right\}\\
\nonumber & = \left( n \lambda_n\right)^{-1}
\left\{\bxi_{\mathfrak{M}_0^c} -
\left[\bgamma_{\mathfrak{M}_0^c}(\hbbeta_1) -
\bgamma_{\mathfrak{M}_0^c}(\bbeta_1)\right] \right \}.
\end{align}
On the event $\mathcal{E}_2$,  the $L_\infty$ norm of the first term
is bounded by $O(n^{-1/2} u_n \lambda_n^{-1}) = o(1)$ by the condition on
$\lambda_n$. It remains to bound the second term of (\ref{116}).

A Taylor expansion of $\bgamma_{\mathfrak{M}_0^c}(\bdelta)$ around $\bbeta_1$ componentwise gives
\begin{equation} \label{117}
\bgamma_{\mathfrak{M}_0^c}(\hbbeta_1) =
\bgamma_{\mathfrak{M}_0^c}(\bbeta_1) +  \bX_2\t
\Sig\left(\btheta_0\right) \bX_1 (\hbbeta_1 - \bbeta_1) + \bw,
\end{equation}
where $\bw = (w_{s+1}, \cdots, w_p)\t$ with $ w_j = \frac{1}{2}
(\hbbeta_1 - \bbeta_1)\t \nabla^2 \gamma_{j}(\bdelta_j)
(\hbbeta_1 - \bbeta_1)$ and $\bdelta_j$ some $s$-vector
lying on the line segment joining $\hbbeta_1$ and $\bbeta_1$. By
(\ref{009}) in Condition \ref{con2} and $\hbbeta_1 \in \mathcal{N}$,
arguing similarly to (\ref{111}), we have
\begin{align} \label{118}
\left\|\bw\right\|_\infty   = O\left[s n^{1 - 2 \gamma} (\log
n)^2\right].
\end{align}
Since $\hbbeta_1$ solves equation $\overline{\bPsi}(\bdelta) =
\bzero$ in (\ref{112}), we have
\begin{equation} \label{119}
\hbbeta_1 - \bbeta_1 =  \left[\bX_1\t \Sig\left(\btheta_0\right)
\bX_1\right]^{-1} \left(\bxi_{\mathfrak{M}_0} - \bet - \br\right).
\end{equation}
It follows from (\ref{007}) and (\ref{008}) in Condition \ref{con2},
(\ref{109}), (\ref{111}), and (\ref{116})--(\ref{119}) that
\begin{align*}
&\left\|\bz\right\|_\infty \leq o(1) + \left( n \lambda_n\right)^{-1} \left\|\bgamma_{\mathfrak{M}_0^c}(\hbbeta_1) - \bgamma_{\mathfrak{M}_0^c}(\bbeta_1)\right\|_\infty \\
& \quad \leq o(1) + \left( n \lambda_n\right)^{-1} \left\|\bX_2\t \Sig\left(\btheta_0\right) \bX_1 \left[\bX_1\t \Sig\left(\btheta_0\right) \bX_1\right]^{-1}\right\|_\infty\\
& \quad \quad \cdot \left(\left\|\bxi_{\mathfrak{M}_0} - \bet\right\|_\infty + \left\|\br\right\|_\infty\right) + \left( n \lambda_n\right)^{-1} \left\|\bw\right\|_\infty \\
& \quad \leq o(1) + \left(n \lambda_n\right)^{-1} O\left\{ n^{\alpha_1} \left[\sqrt{n \log n} + s n^{1 - 2 \gamma} (\log n)^2\right] + s n^{1 - 2 \gamma} (\log n)^2\right\}\\
& \quad \quad + \left\|\bX_2\t \Sig\left(\btheta_0\right) \bX_1
\left[\bX_1\t \Sig\left(\btheta_0\right)
\bX_1\right]^{-1}\right\|_\infty \rho'(d_n) .
\end{align*}
The second term is of order $O(\lambda_n^{-1} n^{-\alpha} (\log n)^{2}) = o(1)$ by (\ref{104}).  Using (\ref{008}), we have
\begin{align*}
\|\bz\|_\infty
& \leq C\rho'(0+) + o(1) < \rho'(0+)
\end{align*}
for sufficiently large $n$.

\smallskip

Finally, note that condition (\ref{023}) for sufficiently
large $n$ is guaranteed by ${\lambda}_n \kappa_0 = o(\tau_0)$ in
Condition \ref{con3}. Therefore, by Theorem~\ref{T1}, we have shown
that $\hbbeta = (\hbbeta_1\t, \hbbeta_2\t)\t$ is a strict local maximizer of the non-concave penalized likelihood $Q_n(\bbeta)$ (\ref{004}) with $\|\hbbeta -
\bbeta_0\|_\infty = O(n^{-\gamma} \log n)$ and $\hbbeta_2 = 0$ under
the event $\mathcal{E}_1 \cap \mathcal{E}_2$.  These along with
(\ref{106}) prove parts a) and
b). This completes the proof.

\subsection{Proof of Theorem \ref{T3}} \label{Sec8.7}
We continue to adopt the notation in the proof of Theorem~\ref{T2}.
To prove the conclusions, it suffices to show that under the given
regularity conditions, there exists a strict local maximizer $\hbbeta$ of
the penalized likelihood $Q_n(\bbeta)$ in (\ref{004}) such that 1)
$\hbbeta_2 = \bzero$ with probability tending to 1 as $n \rightarrow
\infty$ (i.e., sparsity), and 2) $\|\hbbeta_1 - \bbeta_1\|_2 =
O_P(\sqrt{s/n})$ (i.e., $\sqrt{s/n}$-consistency).

\smallskip

\textit{Step 1: Consistency in the $s$-dimensional subspace}. We first constrain $Q_n(\bbeta)$ on the $s$-dimensional subspace $\{\bbeta \in \mathbf{R}^p: \bbeta_{\mathfrak{M}_0^c} = \bzero\}$ of $\mathbf{R}^p$. This constrained penalized likelihood is given by
\begin{equation} \label{126}
\overline{Q}_n(\bdelta) = \overline{\ell}_n(\bdelta) -  \sum_{j =
1}^s p_{\lambda_n}(|\delta_j|),
\end{equation}
where $\overline{\ell}_n(\bdelta) =  n^{-1}[\by\t \bX_1 \bdelta -
\bone\t \bb(\bX_1 \bdelta)]$ and $\bdelta = (\delta_1, \cdots,
\delta_s)\t$. We now show that there exists a strict local maximizer
$\hbbeta_1$ of $\overline{Q}_n(\bdelta)$ such that $\|\hbbeta_1 -
\bbeta_1\|_2 = O_P(\sqrt{s/n})$. To this end, we define an event
\[
H_n = \left\{\overline{Q}_n(\bbeta_1) > \max_{\bdelta \in \partial N_\tau} \overline{Q}_n(\bdelta)\right\},
\]
where $\partial N_\tau$ denotes the boundary of the closed set
$N_\tau = \{\bdelta \in \mathbf{R}^s: \|\bdelta - \bbeta_1\|_2 \leq
\sqrt{s/n} \tau\}$ and $\tau \in (0, \infty)$. Clearly, on the event
$H_n$, there exists a local maximizer $\hbbeta_1$ of
$\overline{Q}_n(\bdelta)$ in $N_\tau$.  Thus, we need only to show that $P(H_n)$
is close to 1 as $n \rightarrow \infty$ when $\tau$ is large.  To this end, we need to analyze
the function $\overline{Q}_n$ on the boundary $\partial N_\tau$.

Let $n$ be sufficiently large such that $\sqrt{s/n} \tau \leq d_n$ since $d_n \gg \sqrt{s/n}$ by Condition \ref{con5}. It
is easy to see that $\bdelta = (\delta_1, \cdots, \delta_s)\t \in
N_\tau$ entails $\sgn(\bdelta) = \sgn(\bbeta_1)$, $\|\bdelta -
\bbeta_1\|_\infty \leq d_n$, and $\min_j |\delta_j| \geq d_n$. By
Taylor's theorem, we have for any $\bdelta \in N_\tau$,
\begin{equation} \label{127}
\overline{Q}_n(\bdelta) - \overline{Q}_n(\bbeta_1) = (\bdelta - \bbeta_1)\t \bv - \frac{1}{2} (\bdelta - \bbeta_1)\t \bD (\bdelta - \bbeta_1),
\end{equation}
where $\bv =  n^{-1} \bX_1\t [\by - \bmu(\btheta_0)] -
\bar{p}_{\lambda_n}(\bbeta_1)$, $\btheta_0 = \bX \bbeta_0 = \bX_1
\bbeta_1$,
\[ \bD = n^{-1} \bX_1\t \Sig\left(\btheta^*\right) \bX_1 +
   \diag\left\{p''_{\lambda_n}(|\bbeta^*|)\right\}, \]
$\btheta^* = \bX_1 \bbeta^*$, and $\bbeta^*$ lies on the line
segment joining $\bdelta$ and $\bbeta_1$. More generally, when the second derivative of the penalty function $p_\lambda$ does not necessarily exist, it is easy to show that the second part of the matrix $\bD$ can be replaced by a diagonal matrix with maximum absolute element bounded by $\lambda_n \kappa_0$. Recall that
\[
\mathcal{N}_0 = \{\bb \in \mathbf{R}^s: \|\bb - \bbeta_1\|_\infty \leq d_n\}
\]
and $\kappa_0 = \max_{\bb \in \mathcal{N}_0} \kappa(\rho; \bb)$,
where $\kappa(\rho; \bb)$ is given by (\ref{016}). For any $\bdelta
\in \partial N_\tau$, we have $\|\bdelta - \bbeta_1\|_2 = \sqrt{s/n}
\tau$ and $\bbeta^* \in \mathcal{N}_0$. Then for sufficiently large
$n$, by (\ref{129}) and $\lambda_n \kappa_0 = o(1)$ in Conditions
\ref{con4} and \ref{con5} we have
\[
\lambda_{\min}(\bD) \geq  c  -  \lambda_n \kappa_0 \geq \frac{c}{2}.
\]
Thus by (\ref{127}), we have
\[
\max_{\bdelta \in \partial N_\tau} \overline{Q}_n(\bdelta) -
\overline{Q}_n(\bbeta_1) \leq \sqrt{s/n} \tau
\left(\left\|\bv\right\|_2 - c \sqrt{s/n} \tau/4\right),
\]
which along with Markov's inequality entails that
\[
P(H_n) \geq P\left(\left\|\bv\right\|_2^2 < \frac{c^2 s \tau^2}{16 n}
\right) \geq 1 - \frac{16  n E \left\|\bv\right\|_2^2}{c^2 s
\tau^2}.
\]

It follows from $E \by = \bmu(\btheta_0)$, $\cov(\by) = \aphi
\Sig(\btheta_0)$, and Conditions \ref{con4} and \ref{con5} that
\begin{align*}
E \left\|\bv\right\|_2^2 & = n^{-2} E \left\|\bX_1\t \left[\by - \bmu\left(\btheta_0\right)\right]\right\|_2^2 + \left\| \bar{p}_{\lambda_n}(\bbeta_1)\right\|_2^2 \\
& \leq  n^{-2} \aphi \tr\left[\bX_1\t \Sig(\btheta_0) \bX_1\right] +
s  p_{\lambda_n}'(d_n)^2 = O(s n^{-1}),
\end{align*}
since $p_{\lambda_n}'(t)$ is decreasing in $t \in [0, \infty)$.
Hence, we have
\[ P(H_n) \geq 1 - O(\tau^{-2}). \]
This proves
$\|\hbbeta_1 - \bbeta_1\|_2 = O_P(\sqrt{s/ n})$.

\smallskip

\textit{Step 2: Sparsity}. Let $\hbbeta \in \mathbf{R}^p$ with $\hbbeta_{\mathfrak{M}_0} = \hbbeta_1 \in N_\tau \subset \mathcal{N}_0 $ a strict local maximizer of $\overline{Q}_n(\bdelta)$ and $\hbbeta_2 = \hbbeta_{\mathfrak{M}_0^c} = \bzero$, and $\hbtheta = \bX \hbbeta$. It remains to prove that the vector $\hbbeta$ is indeed a strict local maximizer of $Q_n(\bbeta)$ on the space $\mathbf{R}^p$. From the proof of Theorem \ref{T1}, we see that it suffices to check condition (\ref{022}). The idea is the same as that in Step 2 of the proof of Theorem \ref{T2}. Let $\bxi = (\xi_1, \cdots, \xi_p)\t = \bX\t \by - \bX\t \bmu(\btheta_0)$ and consider the event
\[ \mathcal{E}_2 = \left\{\left\|\bxi_{\mathfrak{M}_0^c}\right\|_\infty \leq u_n \sqrt{n}\right\}, \]
where $u_n = c_1^{-1/2} n^{\alpha/2} \sqrt{\log n}$. We have shown
in the proof of Theorem \ref{T2} that
\begin{equation} \label{128}
P(\mathcal{E}_2) \geq 1 - (p - s) \varphi(u_n) \geq 1 - 2 p e^{-c_1
u_n^2} \to 1,
\end{equation}
since $\log p = O(n^\alpha)$. It follows from (\ref{130}) and
(\ref{131}) in Condition \ref{con4}, (\ref{116}), (\ref{117}) that
\begin{align*}
\left\|\bz\right\|_\infty & \leq \left( n \lambda_n\right)^{-1} \left[\left\|\bxi_{\mathfrak{M}_0^c}\right\|_\infty + \left\|\bX_2\t \bmu\left(\hbtheta\right) - \bX_2\t \bmu\left(\btheta_0\right)\right\|_\infty\right] \\
& = o(1) +  \left( n \lambda_n\right)^{-1} \left[\left\|\bX_2\t \Sig\left(\btheta_0\right) \bX_1 (\hbbeta_1 - \bbeta_1)\right\|_\infty + \left\|\bw\right\|_\infty\right]\\
& = o(1) + \left( n \lambda_n\right)^{-1}  \left[O(n) \left\|\hbbeta_1 - \bbeta_1\right\|_2 + O(n) \left\|\hbbeta_1 - \bbeta_1\right\|_2^2\right]\\
& = o(1) + O\left(\lambda_n^{-1} \sqrt{s/n} \tau\right) = o(1),
\end{align*}
which shows that inequality (\ref{022}) holds for sufficiently large $n$. This concludes the proof.

\subsection{Proof of Theorem \ref{T4}} \label{Sec8.8}
Clearly by Theorem \ref{T3}, we only need to prove the asymptotic
normality of $\hbbeta_1$. On the event $H_n$ defined in the proof of
Theorem \ref{T3}, it has been shown that $\hbbeta_1 \in N_\tau
\subset \mathcal{N}_0$ is a strict local maximizer of
$\overline{Q}_n(\bdelta)$ and $\hbbeta_2 = \bzero$. It follows
easily that $\nabla \overline{Q}_n(\hbbeta_1) = \bzero$. In view of
(\ref{126}), we have
\[ \nabla \overline{Q}_n(\bdelta) = \nabla
\overline{\ell}_n(\bdelta) - \bar{p}_{\lambda_n}(\bdelta). \] We
expand the first term $\nabla \overline{\ell}_n(\bdelta)$ around
$\bbeta_1$ to the first order componentwise. Then by (\ref{131}) in
Condition \ref{con4} and $\|\hbbeta_1 - \bbeta_1\|_2 =
O_P(\sqrt{s/n})$, we have under the $L_2$ norm,
\begin{align} \label{132}
\bzero & = \nabla \overline{Q}_n(\hbbeta_1) = \nabla \overline{\ell}_n(\bbeta_1)
   - n^{-1} \bX_1\t \Sig\left(\btheta_0\right) \bX_1 (\hbbeta_1 - \bbeta_1) \\
\nonumber
& \quad + O(1) \left\|\hbbeta_1 - \bbeta_1\right\|_2^2 \sqrt{s} -  \bar{p}_{\lambda_n}(\hbbeta_1) \\
\nonumber
& =  n^{-1} \bX_1\t \left[\by - \bmu(\btheta_0)\right] -  n^{-1} \bX_1\t \Sig\left(\btheta_0\right) \bX_1 (\hbbeta_1 - \bbeta_1) \\
\nonumber & \quad -  \bar{p}_{\lambda_n}(\hbbeta_1) + O_P(s^{3/2} n^{-1}).
\end{align}
It follows from $\hbbeta_1 \in \mathcal{N}_0$, and
$p_{\lambda_n}'(d_n) = o(s^{-1/2} n^{-1/2})$ in Condition \ref{con6}
that
\begin{equation} \label{135}
\left\|\bar{p}_{\lambda_n}(\hbbeta_1)\right\|_2 \leq \sqrt{s}
p_{\lambda_n}'(d_n) = o_P(1/\sqrt{n}),
\end{equation}
due to the monotonicity of $p_{\lambda_n}'(t)$. Combing (\ref{132})
and (\ref{135}) gives
\[
\bX_1\t \Sig\left(\btheta_0\right) \bX_1 (\hbbeta_1 - \bbeta_1) =
\bX_1\t \left[\by - \bmu(\btheta_0)\right] + o_P(\sqrt{n}),
\]
since $s = o(n^{1/3})$. This along with the first part of (\ref{129})
in Condition \ref{con4} entails
\begin{equation} \label{136}
\bB_n^{1/2} \left(\hbbeta_1 - \bbeta_1\right)  = \bB_n^{-1/2}\bX_1\t
\left[\by - \bmu(\btheta_0)\right] + o_P(1),
\end{equation}
where $\bB_n = \bX_1\t \Sig(\btheta_0) \bX_1$ and the small order term is understood under the $L_2$ norm.

We are now ready to show the asymptotic normality of $\hbbeta_1$.
Let $\bA_n \bA_n\t \rightarrow \bG$, where $\bA_n$ is a $q \times s$
matrix and $\bG$ is a symmetric positive definite matrix. It follows
from (\ref{136}) that
\[ \bA_n \bB_n^{1/2} \left(\hbbeta_1 - \bbeta_1\right) = \bu_n + o_P(1), \]
where $\bu_n =  \bA_n \bB_n^{-1/2} \bX_1\t [\by - \bmu(\btheta_0)]$.
Thus by Slutsky's lemma, to show that
\[ \bA_n \bB_n^{1/2} \left(\hbbeta_1 - \bbeta_1\right) \toD N(\bzero,  \aphi \bG), \]
it suffices to prove $\bu_n \toD N(\bzero,  \aphi \bG)$. For any
 unit vector $\ba \in \mathbf{R}^q$, we consider the
asymptotic distribution of the linear combination
\[ v_n = \ba\t \bu_n =  \ba\t \bA_n \bB_n^{-1/2} \bX_1\t \left[\by - \bmu(\btheta_0)\right] = \sum_{i = 1}^n \xi_i, \]
where $\xi_i = \ba\t \bA_n \bB_n^{-1/2} \bz_i [y_i - b'(\theta_{0,
i})]$ and $\bX_1 = (\bz_1, \cdots, \bz_n)\t$. Clearly $\xi_i$'s are
independent and have mean 0, and
\begin{align*}
\sum_{i = 1}^n \var(\xi_i) & =  \ba\t \bA_n \bB_n^{-1/2} \aphi \left[\bX_1\t \Sig\left(\btheta_0\right) \bX_1\right] \bB_n^{-1/2} \bA_n\t \ba \\
& =  \aphi \ba\t \bA_n \bA_n\t \ba \longrightarrow  \aphi \ba\t \bG
\ba
\end{align*}
as $n \rightarrow \infty$. By Condition \ref{con6} and the Cauchy-Schwarz inequality, we have
\begin{align*}
\sum_{i = 1}^n E \left|\xi_i\right|^3 & =  \sum_{i = 1}^n \left|\ba\t \bA_n \bB_n^{-1/2} \bz_i\right|^3 E \left|y_i - b'\left(\theta_{0, i}\right)\right|^3 \\
& = O(1) \sum_{i = 1}^n \left|\ba\t \bA_n \bB_n^{-1/2} \bz_i\right|^3 \\
& \leq O(1) \sum_{i = 1}^n \left\|\ba\t \bA_n\right\|_2^3 \left\|\bB_n^{-1/2} \bz_i\right\|_2^3 \\
& = O(1) \sum_{i = 1}^n \left(\bz_i\t \bB_n^{-1} \bz_i\right)^{3/2} = o(1).
\end{align*}
Therefore an application of Lyapunov's theorem yields
\[ \ba\t \bu_n = \sum_{i = 1}^n \xi_i \toD N(0, \aphi \ba\t \bG \ba). \]
Since this asymptotic normality holds for any unit vector $\ba \in
\mathbf{R}^q$, we conclude that $\bu_n \toD N(\bzero, \aphi \bG)$,
which completes the proof.

\appendix
\section{Appendix} \label{SecA}

\subsection{Three commonly used GLMs} \label{SecA.1}
In this section we give the formulas used in the ICA algorithm for three
commonly used GLMs: linear regression model, logistic regression
model, and Poisson regression model.

\smallskip

\textit{Linear regression}. For this model, $b(\theta) = \frac{1}{2} \theta^2$, $\theta \in \mathbf{R}$ and $\aphi = \sigma^2$. The penalized likelihood $Q_n(\bbeta)$ in (\ref{004}) can be written as
\begin{equation} \label{121}
Q_n(\bbeta) = -\left\{(2 n)^{-1} \left\|\by - \bX \bbeta\right\|_2^2
+  \sum_{j = 1}^p
p_\lambda\left(\left|\beta_j\right|\right) \right\},
\end{equation}
where $\bbeta = (\beta_1, \cdots, \beta_p)\t$.  Thus maximizing
$Q_n(\bbeta)$ becomes the penalized least squares problem. In Step 3
of ICA, we have $\widetilde{Q}_n(\beta_j; \hbbeta^{\lambda_k}, j) =
Q_n(\bbeta)$, where the subvector of $\bbeta$ with components in
$\{1, \cdots, p\} \setminus \{j\}$ is identical to that of
$\hbbeta^{\lambda_k}$.

\smallskip

\textit{Logistic regression}. For this model, $b(\theta) = \log(1 +
e^\theta)$, $\theta \in \mathbf{R}$ and $\aphi = 1$. In Step 3 of
ICA, by (\ref{123}) we have
\begin{align} \label{122}
\widetilde{Q}_n&(\beta_j; \hbbeta^{\lambda_k}, j) = \ell_n(\hbbeta^{\lambda_k}) + n^{-1} \left\{\bx_j\t \left[\by - \bmu\left(\bX \hbbeta^{\lambda_k}\right)\right]\right\} \left(\beta_j - \hbeta^{\lambda_k}_j\right) \\
\nonumber
& \quad -\frac{1}{2 n} \left[\bx_j\t \Sig\left(\bX \hbbeta^{\lambda_k}\right) \bx_j\right] \left(\beta_j - \hbeta^{\lambda_k}_j\right)^2 - \sum_{j = 1}^p p_{\lambda_k}\left(\left|\beta_j\right|\right),
\end{align}
where the subvector of $\bbeta$ with components in $\{1, \cdots, p\} \setminus \{j\}$ is identical to that of $\hbbeta^{\lambda_k}$, $\hbbeta^{\lambda_k} = (\hbeta^{\lambda_k}_1, \cdots, \hbeta^{\lambda_k}_p)\t$, $\bX = (\bx_1, \cdots, \bx_p)$, $\bmu(\bX \hbbeta^{\lambda_k}) = \left(\frac{e^{\theta_1}}{1 + e^{\theta_1}},
\cdots, \frac{e^{\theta_n}}{1 + e^{\theta_n}}\right)\t$, and
\[ \Sig(\bX \hbbeta^{\lambda_k}) = \diag\left\{\frac{e^{\theta_1}}{\left(1 + e^{\theta_1}\right)^2}, \cdots,
  \frac{e^{\theta_n}}{\left(1 + e^{\theta_n}\right)^2}\right\} \]
with $(\theta_1, \cdots, \theta_n)\t = \bX \hbbeta^{\lambda_k}$.

\smallskip

\textit{Poisson regression}. For this model, $b(\theta) = e^\theta$,
$\theta \in \mathbf{R}$ and $\aphi = 1$. In Step 3 of ICA,
$\widetilde{Q}_n(\beta_j; \hbbeta^{\lambda_k}, j)$ has the same
expression as in (\ref{122}) with
\[
\bmu(\bX \hbbeta^{\lambda_k}) = \left(e^{\theta_1}, \cdots,
e^{\theta_n}\right)\t \quad \text{and} \quad \Sig(\bX \hbbeta^{\lambda_k}) = \diag\left\{e^{\theta_1}, \cdots, e^{\theta_n}\right\},
\]
where $(\theta_1, \cdots, \theta_n)\t = \bX \hbbeta^{\lambda_k}$.

\subsection{SCAD penalized least squares solution} \label{SecA.2}
Consider the univariate SCAD penalized least squares problem
\begin{equation} \label{124}
\min\limits_{\beta \in \mathbf{R}} \left\{2^{-1} \left(z - \beta\right)^2 + \Lambda p_\lambda\left(\left|\beta\right|\right)\right\},
\end{equation}
where $z \in \mathbf{R}$, $\Lambda \in (0, \infty)$, and $p_\lambda$ is the SCAD penalty given by (\ref{003}). The solution when $\Lambda = 1$ was given by Fan (1997). We denote by $R(\beta)$ the objective function and $\hbeta(z)$ the minimizer of problem (\ref{124}). Clearly $\hbeta(z)$ equals 0 or solves the gradient equation
\begin{equation} \label{125}
g(\beta) \equiv \nabla_\beta \left\{2^{-1} \left(z - \beta\right)^2 + \Lambda p_\lambda\left(\left|\beta\right|\right)\right\} = \beta - z + \sgn(\beta) \Lambda p_\lambda'(|\beta|) = 0.
\end{equation}
It is easy to show that $\hbeta(z) = \sgn(z) |\hbeta(z)|$ and $|\hbeta(z)| \leq |z|$, i.e., $\hbeta(z)$ is between 0 and $z$. Let $z_0 = \sgn(z) (|z| - \Lambda \lambda)_+$.

\smallskip

1) If $|z| \leq \lambda$, we can easily show that $\hbeta(z) = z_0$.

\smallskip

2) Let $\lambda < |z| \leq a \lambda$. Note that $g$ defined in (\ref{125}) is piecewise linear between 0 and $z$, and $g(0) = \sgn(z)[-|z| + \Lambda \lambda]$, $g(\sgn(z) \lambda) = \sgn(z)[-|z| + (\Lambda + 1) \lambda]$, $g(z) = \sgn(z) \Lambda p_\lambda'(|z|)]$. Thus it is easy to see that if $|z| \leq (\Lambda + 1) \lambda$, we have $\hbeta(z) = z_0$, and if $|z| > (\Lambda + 1) \lambda$, we have
\[ \hbeta(z) = \sgn(z) \frac{|z| - \Lambda \lambda (a - 1)^{-1} a}{1 - (a - 1)^{-1} \Lambda}. \]

3) Let $|z| > a \lambda$. The same argument as in 2) shows that when $|z| \leq (\Lambda + 1) \lambda$, we have $\hbeta(z) = z_0$ if $R(z_0) \leq R(z)$ and $\hbeta(z) = z$ otherwise. When $|z| > (\Lambda + 1) \lambda$, we have $\hbeta(z) = z$.

\begin{singlespace}

\end{singlespace}

\end{document}